\documentclass[a4paper]{article}

\usepackage{amsmath, amssymb, mathrsfs}
\usepackage{graphicx, subfigure, epsfig, tikz}
\usepackage{booktabs}
\usepackage{hyperref}

\graphicspath{{Figure/}}

\newtheorem{problem}{Problem}
\newtheorem{remark}{Remark}
\newtheorem{example}{Example}

\begin{document}

\title{A reference ball based iterative algorithm for imaging acoustic obstacle from phaseless far-field data}

\author{
Heping Dong\thanks{School of Mathematics, Jilin University, Changchun, P. R. China. {\it dhp@jlu.edu.cn}}, 
Deyue Zhang\thanks{School of Mathematics, Jilin University, Changchun, P. R. China. {\it dyzhang@jlu.edu.cn} (Corresponding author)}\ \ and 
Yukun Guo\thanks{Department of Mathematics, Harbin Institute of Technology, Harbin, P. R. China. {\it ykguo@hit.edu.cn}}}

\maketitle

\begin{abstract}
In this paper, we consider the inverse problem of determining the location and the shape of a sound-soft obstacle from the modulus of the far-field data for a single incident plane wave. By adding a reference ball artificially to the inverse scattering system, we propose a system of nonlinear integral equations based iterative scheme to reconstruct both the location and the shape of the obstacle. The reference ball technique causes few extra computational costs, but breaks the translation invariance and brings information about the location of the obstacle. Several validating numerical examples are provided to illustrate the effectiveness and robustness of the proposed inversion algorithm.
\end{abstract}



\section{Introduction}

It is well known that inverse scattering problems (ISP) play a significantly crucial role in a vast variety of realistic applications such as radar sensing, ultrasound tomography, biomedical imaging, geophysical exploration and noninvasive detecting. Therefore, effective and efficient numerical inversion approaches have been extensively and intensively studied in the recent decades \cite{DR-shu2}.

In typical inverse scattering problems, reconstruction of the geometrical information is based on the knowledge of full measured data (both the intensity and phase information are collected). However, in a number of practical scenarios, to measure the full data is extremely difficult or might even be unavailable. Thus, phaseless inverse scattering problems arise naturally and attract a great attention from mathematical and numerical point of view. Without the phase information, the theoretical justifications of the uniqueness issue and the development of inversion algorithms are more challenging than the phased inverse scattering problem.

Various numerical methods have been proposed for solving phaseless inverse acoustic obstacle scattering problems. Kress and Rundell \cite{RW1997} investigate the phaseless inverse obstacle scattering and propose a Newton method for imaging a two-dimensional sound-soft obstacle from the modulus of the far-field data with only one incoming direction. In particular, it is pointed out that the location of the obstacle cannot be reconstructed since the modulus of the far-field pattern has translation invariance. This means that the solution of the inverse problem is not unique. If the location of the obstacle is known, then a uniqueness result is presented in \cite{LZ2010} for a special case, that is, the sound-soft ball can be determined by phaseless far-field pattern when the product of wavenumber and the radius of the ball is less than a constant. In \cite{Ivanyshyn2007}, a nonlinear integral equation method is investigated for the two-dimensional shape reconstruction from a given incident field and the modulus of the far-field data. This method involves the full linearization of a pair of integral equations, i.e., field equation and phaseless data equation, with respect to both the boundary parameterization and the density. Then, the nonlinear integral equation method is extended to the three-dimensional shape reconstruction in \cite{OR2010}. The problem is divided into two subproblems including the shape reconstruction from the modulus of the far-field data and the location identification based on the invariance relation using only few far field measurements in the backscattering direction. In addition, fundamental solution method \cite{KarageorghisAPNUM} and a hybrid method \cite{Lee2016} are proposed to detect the shape of a sound-soft obstacle by using of the modulus of the far-field data for one incident field. 

Besides the aforementioned phaseless inverse acoustic obstacle scattering, there exist other types of phaseless inverse scattering problems as well as the relevant numerical methods. Based on the physical optics approximation and the local maximum behavior of the backscattering far-field pattern, a numerical method is developed in \cite{Li2017} to reconstruct an electromagnetic polyhedral PEC obstacle from a few phaseless backscattering measurements. In \cite{CH2016}, a direct imaging method based on reverse time migration is proposed for reconstructing extended obstacles by phaseless electromagnetic scattering data.  A continuation method \cite{BLL2013} is developed to reconstruct the periodic grating profile from phaseless scattered data. In \cite{Bao2016}, a recursive linearization algorithm is proposed to image the multi-scale rough surface with tapered incident waves from measurements of the multi-frequency phaseless scattered data. In \cite{ChengJin2017}, an iterative method based on Rayleigh expansion is developed to solve the inverse diffraction grating problem with super-resolution by using phase or phaseless near-field data. For a recent work on the Fourier method for solving phaseless inverse source problem, we refer to \cite{ZGLL18}. We also refer to \cite{Ammari2016, Kli2017, KR2016a, KR2016b} for phaseless inverse medium scattering problems.

A major difficulty in solving the phaseless inverse obstacle scattering problem is to determine the location of the obstacle from phaseless farfield data, since the modulus of the far-field pattern is invariant under translations. Recently, a recursive Newton-type iteration method is developed in \cite{ZhangBo2017} to recover both the location and the shape of the obstacle simultaneously from multi-frequency phaseless far-field data. In this numerical method, some sets of superpositions of two plane waves with different directions are used as the incident fields to break the translation invariance property of the phaseless far-field pattern.

In this paper, we consider the inverse problem of determining the location and the shape of a sound-soft obstacle from the modulus of the far-field data for a single incident plane wave. In a recent work \cite{GDM2018}, the nonlinear integral equation method proposed by Johansson and Sleeman \cite{TB2007} is extended to reconstruct the shape of a sound-soft crack by using phaseless far-field data for one incident plane wave. Motivated by \cite{GDM2018, TB2007} and the reference ball technique in \cite{LiJingzhi2009, ZG18}, we first introduce an artificial reference ball to the inverse scattering system, and propose an iterative scheme which involves a system of nonlinear and ill-posed integral equations to reconstruct both the location and the shape of the obstacle. Since the location of the reference ball is known and fixed, it has the capability of calibrating the scattering system so that the translation invariance does not occur. As a result, the location information of the obstacle could be recovered with negligible additional computational costs.

To our best knowledge, this is the first attempt in the literature towards reconstructing both the location and the shape of the unknown obstacle by using the intensity-only far-field data due to a single incident plane wave. The novelty of this work lies in the incorporation of the reference ball technique into the iterative scheme. Hence, our approach exhibits the following salient features: First, rather than the superposition of different incident waves, only a single incident plane wave with a fixed wavenumber and a fixed incoming direction is needed. Second, the iteration procedure is very fast and easy to implement because the Fr\'{e}chet derivatives can be explicitly formulated and thus do not require the forward solver. Third, the location and profile information of the obstacle could be simultaneously reconstructed. Finally, the iterative scheme is robust in the sense that it is insensitive to the parameters such as the initial guess, the location and size of the reference ball, as well as the measurement noise.

The rest of this paper is organized as follows. In the next section, we formulate the phaseless inverse obstacle scattering problem in conjunction with the reference ball. In Section 3, we derive a system of boundary integral equations with a reference ball, and an iterative scheme is presented to solve the boundary integral equations. In Section 4, several numerical examples will be presented, including the numerical implementation details. Finally, some concluding remarks are summarized in Section 5.


\section{Problem formulation}

Let $D\subset\mathbb{R}^2$ be an open bounded domain with $\mathcal{C}^2$ boundary and the positive constant $\kappa$ be the wavenumber. Given an incident plane wave $u^i=\mathrm{e}^{\mathrm{i}\kappa x\cdot d}$ with the incoming direction $d$, the forward/direct obstacle scattering problem is to find the total field $u$ satisfying the Helmholtz equation
\begin{equation}\label{Helmequ}
\Delta u + \kappa^2u = 0, \quad\text{in}~\mathbb{R}^2\backslash\overline{D}, 
\end{equation}
and the Dirichlet boundary condition
\begin{equation}\label{DirichletBC}
u=0, \quad \text{on}\ \partial D. 
\end{equation}
The total field $u=u^i+u^s$ is given as the sum of the known incident wave $u^i$ and the unknown scattered wave $u^s$ which is required to fulfill the Sommerfeld radiation condition
\begin{equation}\label{SRC}
\lim_{r=|x|\to\infty}r^{1/2}\left(\frac{\partial u^s}{\partial r}-\mathrm{i}\kappa u^s\right)=0. 
\end{equation}

Further, the scattered field $u^s$ has the following asymptotic behavior \cite{DR-shu2}
$$
u^s(x) = \frac{\mathrm{e}^{\mathrm{i}\kappa|x|}}{\sqrt{|x|}}\left\{u^\infty(\hat{x}) +
\mathcal{O}\left(\frac{1}{|x|}\right)\right\},~~~\text{as}~|x|\rightarrow\infty,
$$
uniformly for all direction $\hat{x}:=x/|x|$, where the complex-valued function $u^\infty(\hat{x})$ defined on the unit circle $\Omega$ is known as the far field pattern or scattering amplitude. Then, the phaseless inverse scattering problem is stated as follows:

\begin{problem}[Phaseless ISP]\label{problem_1}
	Given an incident plane wave $u^i$ for a single wavenumber and a single incident direction $d$, together with the corresponding phaseless far-field data $|u_D^\infty(\hat x)|, \  \hat x\in\Omega,$ due to the unknown obstacle $D$, determine the location and shape of $\partial D$.
\end{problem}

It has been pointed out in \cite{RW1997} that Problem \ref{problem_1} does not admit a unique solution, namely, the location of the obstacle cannot be reconstructed since the modulus of the far-field pattern has translation invariance. Specifically, for the shifted domain $D^h:=\{x+h: x\in D\}$ with a fixed vector $h\in\mathbb{R}^2$, the far-field pattern $u_{D^h}^\infty$ satisfies the relation
$$
u_{D^h}^\infty(\hat{x},d)=\mathrm{e}^{\mathrm{i}\kappa h\cdot(d-\hat{x})}u_{D}^\infty(\hat{x},d).
$$
Moreover, the ambiguity induced by this translation invariance relation cannot be remedied by using a finite number of incident waves with different wavenumbers or different incident directions, see \cite{RW1997}. Our goal in this paper is to overcome this difficulty via the reference ball based iterative scheme. To this end, we reformulate Problem \ref{problem_1} as follows:

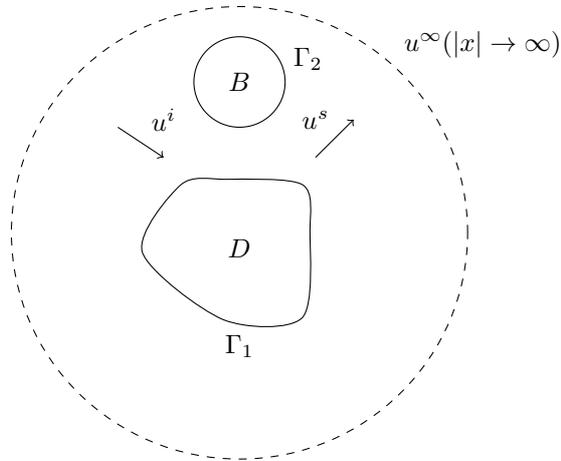
\begin{figure}
	\centering
	\begin{tikzpicture}
	\pgfmathsetseed{8}
	\draw plot [smooth cycle, samples=7, domain={1:8}] (\x*360/8+5*rnd:0.5cm+1cm*rnd) node at (0,-0.2) {$D$}; 
	\draw node at (0, -1.5) {$\Gamma_1$};
	\draw (0, 2) circle (0.6cm) node at (0,2) {$B$}; 
	\draw node at (0.9, 2.3) {$\Gamma_2$};
	\draw [->] (-1.6,1.4)--(-1,1) node at (-1,1.5) {$u^i$}; 
	\draw [->] (1,1)--(1.5, 1.5) node at (1,1.5) {$u^s$}; 
	\draw [dashed] (0, 0) circle (3cm) node at (3.2, 2.5) {$u^\infty(|x|\to \infty)$}; 
	\end{tikzpicture}
	\caption{An illustration of the reference ball technique.} \label{fig:illustration}
\end{figure}

\begin{problem}[Phaseless ISP with a reference ball]\label{problem_2}\ Let $B\subset\mathbb{R}^2$ be the artificially added sound-soft ball such that $D\cap B=\emptyset$. Given an incident plane wave $u^i$ for a single wavenumber and a single incident direction $d$, together with the corresponding phaseless far-field data $|u_{D\cup B}^\infty(\hat x)|, \hat x\in\Omega,$ due to the scatterer $D\cup B$, determine the location and shape of $\partial D$.
\end{problem}

The geometry setting of Problem \ref{problem_2} is illustrated in Fig. \ref{fig:illustration}. For brevity, we denote by $\Gamma_1:=\partial D$ and $\Gamma_2:=\partial B$ in what follows. In the next section, we will introduce a system of nonlinear integral equations based iterative scheme for solving Problem \ref{problem_2}. After that, several numerical experiments will be conducted to demonstrate the feasibility of breaking the translation invariance using the proposed method.

\section{The inversion scheme}

\subsection{Nonlinear integral equations}

We begin this section with establishing a system of nonlinear integral equations for the inverse scattering problem. Denote the fundamental solution to the Helmholtz equation by
$$
\Phi(x,y)=\frac{\mathrm{i}}{4}H_0^{(1)}(\kappa|x-y|), \quad x\neq y,
$$
where $H_0^{(1)}$ denotes the zero-order Hankel function of the first kind. By applying the Huygens' principle \cite[Theorem 3.14]{DR-shu2} to the scattering system $D\cup B$, we obtain
\begin{equation}
u(x)= u^i(x)-\sum\limits_{j=1,2}\int_{\Gamma_j}\frac{\partial u}{\partial\nu}(y)\Phi(x,y)\,\mathrm{d}s(y), \quad x\in \mathbb{R}^2\backslash\overline{D\cup B}. \label{Huygens1}
\end{equation}
Then the far-field pattern of the scattered field $u^s$ is given by
\begin{equation}
u^\infty(\hat{x})=-\gamma\sum\limits_{j=1,2}\int_{\Gamma_1}\frac{\partial u}{\partial
	\nu}(y)\mathrm{e}^{-\mathrm{i}\kappa\hat{x}\cdot y}\,\mathrm{d}s(y),\quad\hat{x}\in\Omega. \label{Huygens2}
\end{equation}
where $\gamma=\mathrm{e}^{\mathrm{i}\pi/4}/\sqrt{8\kappa\pi}$ and $\nu$ denotes the unit outward normal. We define the single-layer operator
$$
(S_{jl} g)(x):=\int_{\Gamma_j}\Phi(x,y) g(y)\,\mathrm{d}s(y), \quad x\in\Gamma_l, \quad j,l=1,2
$$
and the corresponding far-field operator
$$
(S^\infty_{j}g)(\hat{x}):=-\gamma\int_{\Gamma_j}\mathrm{e}^{-\mathrm{i}\kappa\hat{x}\cdot
	y}g(y)\,\mathrm{d}s(y), \quad \hat{x}\in\Omega,\quad j=1,2.
$$

Let $x\in \mathbb{R}^2\backslash\overline{D\cup B}$ tend to the boundary $\Gamma_1$ and $\Gamma_2$, respectively, and in terms of \eqref{Huygens1}, the Dirichlet boundary condition \eqref{DirichletBC} and the continuity of the single-layer potential, one can readily deduce the following field equations 
\begin{align}
S_{11}g_1+S_{21}g_2 & =u^i\quad \textrm{on}\ \Gamma_1, \label{MHuygens1} \\
S_{12}g_1+S_{22}g_2 & =u^i\quad \textrm{on}\ \Gamma_2, \label{MHuygens3}
\end{align}
with respect to the densities $g_j=\partial u/{\partial\nu|_{\Gamma_j}}, j=1,2$. 

On the other hand, let $|x|\to\infty$, then the linearity of the forward scattering problem and equation \eqref{Huygens2} lead to the phaseless data equation
\begin{equation}\label{MHuygens2}
|S^\infty_1 g_1+S^\infty_2 g_2|^2=|u_{D\cup B}^\infty|^2 \quad\textrm{on}\ \Omega. 
\end{equation}


\subsection{The iterative procedure}

We now seek a sequence of approximation to $\Gamma_1$ by solving the field equations \eqref{MHuygens1}--\eqref{MHuygens3} and phaseless data equation \eqref{MHuygens2} in an alternating manner. Given an approximation for the boundary $\Gamma_1$ one can solve \eqref{MHuygens1}--\eqref{MHuygens3} for $g_1$ and $g_2$. Then keeping $g_1$ and $g_2$ fixed, the equation \eqref{MHuygens2} is linearized with respect to $\Gamma_1$ to update the boundary approximation.

For the sake of simplicity, the boundary $\Gamma_1$ is assumed to be a starlike curve with the parametrized form
\begin{equation}\label{obstacle}
\Gamma_1=\{p_1(\hat{x})=c+r(\hat{x})\hat{x}: ~c=(c_1,c_2),\ \hat{x}\in\Omega\}, 
\end{equation}
the boundary $\Gamma_2$ and the unit circle $\Omega$ are parameterized by
$$
\Gamma_2=\{p_2(\hat{x})=b+R\hat{x}: ~b=(b_1,b_2),\ \hat{x}\in\Omega\},
$$
where
$$
\Omega=\{\hat{x}(t)=(\cos t, \sin t): ~0\leq t< 2\pi\}.
$$
Let $p_j(t)$ be the points on $\Gamma_j$, described by $p_1(t):=(c_1,c_2)+r(t)(\cos t, \sin t)$ and $p_2(t):=(b_1,b_2)+R(\cos t, \sin t)$, where $0\leq t<2\pi$. Further, we introduce the parameterized single-layer operators $S_{jl} $ and far field operator $S^\infty_j$ by
\begin{align}\label{A}
A_{jl}(p_l,\psi_j)(t)=\frac{\mathrm{i}}{4}\int_0^{2\pi}H_0^{(1)}(\kappa|p_l(t)-p_j(\tau)|)\psi_j(\tau)\,\mathrm{d}\tau,\quad j=1,2
\end{align}
and
\begin{align}\label{Ainfty}
A^\infty_j(p_j,\psi_j)(t)&=-\gamma\int_0^{2\pi}\mathrm{e}^{-\mathrm{i}\kappa\hat{x}(t)\cdot p_j(\tau)}\psi_j(\tau)\,\mathrm{d}\tau,\quad j=1,2
\end{align}
where we have set $\psi_j(\tau)=G_j(\tau)g_j(p_j(\tau))$. Here, $G_1(\tau)=\big(r^2(\tau)+\big(\frac{dr}{d\tau}(\tau)\big)^2\big)^{1/2}$ and $G_2(\tau)=R$ denote the Jacobian of the transformation. The corresponding right-hand sides
are $w_l(t)=u^i(p_l(t))$ and $w^\infty(t)=u_{D\cup B}^\infty(\hat{x}(t))$. Thus we can obtain the parametrized integral equations \eqref{MHuygens1}--\eqref{MHuygens2} in the form
\begin{align}
A_{11}(p_1,\psi_1)+A_{21}(p_1,\psi_2) & =w_1, \label{PHuygens1}\\
A_{12}(p_2,\psi_1)+A_{22}(p_2,\psi_2) & =w_2, \label{PHuygens3}\\
|A^\infty_1(p_1,\psi_1)+A^\infty_2(p_2,\psi_2)|^2 & =|w^\infty|^2.
\label{PHuygens2}
\end{align}
The linearization of equation \eqref{PHuygens2} with respect to $p_1$ requires the Fr\'{e}chet derivatives of the operator
$A_1^\infty$, that is
\begin{align}\label{FAinfty}
\left(A'^\infty_1[p_1,\psi_1]q\right)(t) = & \mathrm{i}\kappa\gamma\int_0^{2\pi}
\mathrm{e}^{-\mathrm{i}\kappa\hat{x}(t)\cdot p_1(\tau)}\hat{x}(t)\cdot q(\tau)\psi_1(\tau)\,\mathrm{d}\tau \nonumber \\
= & \mathrm{i}\kappa\gamma\int_0^{2\pi}\exp(-\mathrm{i}\kappa(c_1\cos t+c_2\sin t+r(\tau)\cos(t-\tau))) \nonumber \\
&\qquad\quad\cdot(\Delta c_1\cos t+\Delta c_2\sin t +\Delta r(\tau) \cos(t-\tau))\psi_1(\tau)\,\mathrm{d}\tau. 
\end{align}
where the update $q(\tau)=(\Delta c_1, \Delta c_2)+\Delta r(\tau)(\cos\tau,\sin\tau)$. Hence, by using of the product rule, we have
\begin{align*}
&\left(\overline{A^\infty_1(p_1,\psi_1)+A^\infty_2(p_2,\psi_2)}
\big(A^\infty_1(p_1,\psi_1)+A^\infty_2(p_2,\psi_2)\big)\right)'q \\
=&2\Re\left(\overline{A^\infty_1(p_1,\psi_1)+A^\infty_2(p_2,\psi_2)}
A'^\infty_1[p_1,\psi_1]q\right),
\end{align*}
and the linearization of \eqref{PHuygens2} leads to
\begin{equation}\label{LHuygens2}
Bq=f.
\end{equation}
where
\begin{align*}
Bq:= & 2\Re\left(\overline{A^\infty_1(p_1,\psi_1)+A^\infty_2(p_2,\psi_2)}A'^\infty_1[p_1,\psi_1]q\right) \\
f:= & |w^\infty|^2-|A^\infty_1(p_1,\psi_1)+A^\infty_2(p_2,\psi_2)|^2.
\end{align*}

As usual for iterative algorithms, the stopping criteria is necessary to justify the convergence in numerics. Regarding our iterative procedure, we choose the following relative error estimator
\begin{equation}
E_k:=\frac{\||w^\infty|^2-|A^\infty_1(p^{(k)}_1,\psi_1)+A^\infty_2(p_2,\psi_2)|^2\|_{L^2}}
{\||w^\infty|^2\|_{L^2}}\leq\epsilon \label{relativeerror}
\end{equation}
for some sufficiently small parameter $\epsilon>0$ depending on the noise level. Here, $p^{(k)}_1$ is the $k$th approximation of the boundary $\Gamma_1$.

We are now in the position to present the reference ball based iteration algorithm:
\begin{table}[ht]
	\centering
	\begin{tabular}{cp{.8\textwidth}}
		\toprule
		\multicolumn{2}{l}{{\bf Algorithm:}\quad Iterative procedure for phaseless inverse scattering} \\
		\midrule
		{\bf Step 1} & With a mild a priori information of the unknown scatterer $D$, add a suitable reference ball $B$ such that $D\cap B=\emptyset$; \\ 
		{\bf Step 2} & Emanate an incident plane wave with a fixed wavenumber $\kappa>0$ and a fixed incident direction $d\in\Omega$, and then collect the corresponding noisy phaseless far-field data $|u_{D\cup B}^\infty(\hat x)|, \hat x\in\Omega$ for the scatterer $D\cup B$; \\
		{\bf Step 3} & Select an initial star-like curve $\Gamma^{(0)}$ for the boundary $\partial D$ and the error tolerance $\epsilon$. Set $k=0$; \\
		{\bf Step 4} & For the curve $\Gamma^{(k)}$, find the densities $\psi_1$ and $\psi_2$ from \eqref{PHuygens1}--\eqref{PHuygens3}. \\
		{\bf Step 5} & Solve \eqref{LHuygens2} to obtain the updated approximation $\Gamma^{(k+1)}:=\Gamma^{(k)}+q$ and evaluate the error $E_{k+1}$ defined in \eqref{relativeerror}; \\
		{\bf Step 6} & If $E_{k+1}\geq\epsilon$, then set $k=k+1$ and go to Step 4. Otherwise, the current approximation $\Gamma^{(k+1)}$ is served as the final reconstruction of $\partial D$. \\
		\bottomrule
	\end{tabular}
\end{table}

\begin{remark}
	The advantage of introducing the reference ball is that we would be able to obtain the location of the obstacle via this iterative method, since the update information $(\Delta c_1, \Delta c_2)$ about the location of the obstacle is contained in the term $\Re(\overline{A^\infty_2(p_2,\psi_2)}A'^\infty_1[p_1,\psi_1]q)$ of equation \eqref{LHuygens2}.
\end{remark}

\begin{remark}
	For typical nonlinear integral equations based iterative schemes (see, e.g., \cite{Ivanyshyn2007}), the linearization are carried out with respect to both the boundary and the density functions, so the process of solving unknowns should be essentially intertwined. In contrast, we would like to emphasize that our inversion scheme requires the linearization only with respect to the boundary for phaseless data equation, hence the field equations and the phaseless data equation could be completely decoupled and then solved alternatively and separately in each iterative loop. Moreover, the Fr\'{e}chet derivatives in this paper also does not rely on any solution process and can be explicitly given. Therefore, the numerical implementation is now significantly simplified.
\end{remark}
\section{Numerical experiments}
\subsection{Discretization}

In the following, we describe the fully discretizations of \eqref{PHuygens1}--\eqref{PHuygens3} and \eqref{LHuygens2}
respectively. Let $\tau_j^n:=\pi j/n$, $j=0,\cdots,2n-1$ be an equidistant set of quadrature knots. Setting $w_{l,s}^n=w_l(\tau_s^n)$ and $\psi_{l,j}^n=\psi_l(\tau_j^n)$ for $s,j=0,\cdots,2n-1$ and $l=1,2$.
By using the Nystr\"{o}m method \cite{DR-shu2}, we see that the full discretization of \eqref{PHuygens1}--\eqref{PHuygens3} is of the form
\begin{align*}
w_{1,s}^n= & \sum_{j=0}^{2n-1}\left(R_{|s-j|}^{n}K^1_1(\tau_s^n,\tau_j^n)+
\frac{\pi}{n}K^1_2(\tau_s^n,\tau_j^n)\right)\psi_{1,j}^n \\
& +\frac{\pi}{n}\sum_{j=0}^{2n-1}\frac{\mathrm{i}}{4}
H_0^{(1)}(\kappa|p_1(\tau_s^n)-p_2(\tau_j^n)|)\psi_{2,j}^n, 
\end{align*}
\begin{align*}
w_{2,s}^n= & \sum_{j=0}^{2n-1}\big(R_{|s-j|}^{n}K^2_1(\tau_s^n,\tau_j^n)
+\frac{\pi}{n}K^2_2(\tau_s^n,\tau_j^n)\big)\psi_{2,j}^n \\
&+\frac{\pi}{n}\sum_{j=0}^{2n-1}\frac{\mathrm{i}}{4}
H_0^{(1)}(\kappa|p_2(\tau_s^n)-p_1(\tau_j^n)|)\psi_{1,j}^n. 
\end{align*}
where
\begin{align*}
R_j^{n}:=&-\frac{2\pi}{n}\sum_{m=1}^{n-1}\frac{1}{m}\cos\frac{mj\pi}{n}-\frac{(-1)^j\pi}{n^2}, \\
K^l_1(t,\tau)=& -\dfrac{1}{4\pi}J_0(k|p_l(t)-p_l(\tau)|), \quad l=1,2, \\
K^l_2(t,\tau)=& K^l(t,\tau)-K^l_1(t,\tau)\ln\big(4\sin^2\frac{t-\tau}{2}\big), \quad l=1,2,
\end{align*}
and the diagonal term can be deduced as
$$
K^l_2(\tau,\tau)=\left\{\frac{i}{4}-\frac{E_c}{2\pi}-\frac{1}{2\pi}\ln\big(\frac{\kappa}{2}G_l(\tau)\big)\right\},\quad l=1,2
$$
with the Euler constant $E_c=0.57721\cdots$.

We proceed by discussing the discretization of the linearized equation \eqref{LHuygens2} using Newton's method with least squares \cite{Kress2003}. As finite dimensional space for the approximation of the radial function $r$ and its update $\Delta r$ we choose the space of trigonometric polynomials of the form
\begin{equation}\label{updataq}
\Delta r(\tau)=\sum_{m=0}^M\alpha_m\cos{m\tau}+\sum_{m=1}^M\beta_m\sin{m\tau}. 
\end{equation}
where the integer $M>1$ signifies the truncation.

For simplicity, we reformulate the parameterized operators $A_{jl}, A^\infty_j$ by introducing the following definitions
\begin{align*}
M_1(t,\tau)=&-\gamma\exp(-\mathrm{i}\kappa(c_1\cos t+c_2\sin t+r(\tau)\cos(t-\tau))\psi_1(\tau),\\
M_2(t,\tau)=&-\gamma\exp(-\mathrm{i}\kappa(b_1\cos t+b_2\sin t+R\cos(t-\tau))\psi_2(\tau),\\
L_1(t,\tau)=&\mathrm{i}\kappa\gamma\exp(-\mathrm{i}\kappa(c_1\cos t+c_2\sin t+r(\tau)\cos(t-\tau))\cos t~ \psi_1(\tau),\\
L_2(t,\tau)=&\mathrm{i}\kappa\gamma\exp(-\mathrm{i}\kappa(c_1\cos t+c_2\sin t+r(\tau)\cos(t-\tau))\sin t~ \psi_1(\tau),\\
L_3(t,\tau)=&\mathrm{i}\kappa\gamma\exp(-\mathrm{i}\kappa(c_1\cos t+c_2\sin t+r(\tau)\cos(t-\tau))\cos(t-\tau)~\psi_1(\tau),
\end{align*}
Then, by combining \eqref{Ainfty}, \eqref{FAinfty} and \eqref{LHuygens2} together, we obtain the linear system
\begin{align}
f(\tau_s^n)=&(B_1\Delta c_1)(\tau_s^n)+(B_2\Delta c_2)(\tau_s^n)\nonumber \\
&+\sum_{m=0}^M\alpha_m(B_3\cos{m\tau})(\tau_s^n)+\sum_{m=1}^{M}\beta_m(B_3\sin{m\tau})(\tau_s^n)\label{RLHuygens2}
\end{align}
to be solved to determine the real coefficients $\Delta c_1$,
$\Delta c_2$, $\alpha_m$ and $\beta_m$, where
\begin{align*}
(B_i\chi(\tau))(t)=2\Re\left\{\overline{\int_0^{2\pi}
	\Big(M_1(t,\tau)+M_2(t,\tau)\Big)\mathrm{d}\tau}\int_0^{2\pi}L_i(t,\tau)
\chi(\tau)\mathrm{d}\tau\right\}
\end{align*}
with $i=1,2,3$. In general, $2M+1\ll 2n$ and due to the ill-posedness, the overdetermined system \eqref{RLHuygens2} is solved via the Tikhonov regularization, that is, by minimizing the penalized defect
\begin{align}
& \sum_{s=0}^{2n-1}\big|(B_1\Delta c_1)(\tau_s^n)+(B_2\Delta c_2)(\tau_s^n) \nonumber \\
&+\sum_{m=0}^M\alpha_m(B_3\cos{m\tau})(\tau_s^n) 
+\sum_{m=1}^{M}\beta_m(B_3\sin{m\tau})(\tau_s^n)-f(\tau_s^n)\big|^2 \nonumber \\ 
& +\lambda\bigg(|\Delta c_1|^2+|\Delta c_2|^2+2\pi\big[\alpha_0^2+\frac{1}{2}\sum_{m=1}^M(1+m^2)^2(\alpha_m^2+\beta_m^2)\big]\bigg) \label{RLHuygens3}
\end{align}
with a positive regularization parameter $\lambda$ and $H^2$ penalty term.

In view of the trapezoidal rule
\begin{equation}\label{traperule}
\int_0^{2\pi}f(\tau)\,\mathrm{d}\tau\approx\frac{\pi}{n}\sum_{j=0}^{2n-1}f(\tau_j^n),
\end{equation}
one obtains the approximation of $(B_i\chi(\tau))(\tau_s^n)$, that is, for $i=1,2,3$,
\begin{align*}
\widetilde{B}_i^{\chi}(\tau_s^n)\approx2\left(\frac{\pi}{n}\right)^2\sum_{s=0}^{2n-1}\sum_{j=0}^{2n-1}
\Re\left\{\overline{M_1(\tau_s^n,\tau_j^n)+M_2(\tau_s^n,\tau_j^n)}
L_i(\tau_s^n,\tau_j^n)\right\}\chi(\tau_j^n)
\end{align*}
where the function $\chi(\tau)$ is chosen to be the constant $\Delta c_1$, $\Delta c_2$ or the function $\cos(m\tau)$, $\sin(m\tau)$. Then, the approximations are denoted in turn by $\widetilde{B}_1$, $\widetilde{B}_2$,$\widetilde{B}_3^{cos,m}$ and $\widetilde{B}_3^{sin,m}$. Furthermore, we assume that
\begin{align*}
\widetilde{B}=&(\widetilde{B}_1,\widetilde{B}_2,\widetilde{B}_3^{cos, 0},\cdots,\widetilde{B}_3^{cos,M},\widetilde{B}_3^{sin,1},
\cdots,\widetilde{B}_3^{sin,M})_{(2n)\times(2M+3)}\\
\xi=&(\Delta c_1, \Delta c_2, \alpha_0,\cdots, \alpha_M,\beta_1,\cdots, \beta_M)^\top \\
\widetilde{I}= & \mathrm{diag}\{1, 1, 2\pi, \pi(1+1^2)^2, \cdots, \pi(1+M^2)^2,
\pi(1+1^2)^2, \cdots, \pi(1+M^2)^2\} \\
\widetilde{f}=&(f(\tau_0^n),\cdots, f(\tau_{2n-1}^n))^\top.
\end{align*}
Thus, the minimizer in \eqref{RLHuygens3} is equivalent to the unique solution of
\begin{align}\label{EqualRLHuygens3}
\lambda \widetilde{I}\xi+\widetilde{B}^*\widetilde{B}\xi=\widetilde{B}^*\widetilde{f}.
\end{align}

\subsection{Numerical examples}

\begin{figure}
	\centering 
	\subfigure[5th iteration]{\includegraphics[width=0.32\textwidth]{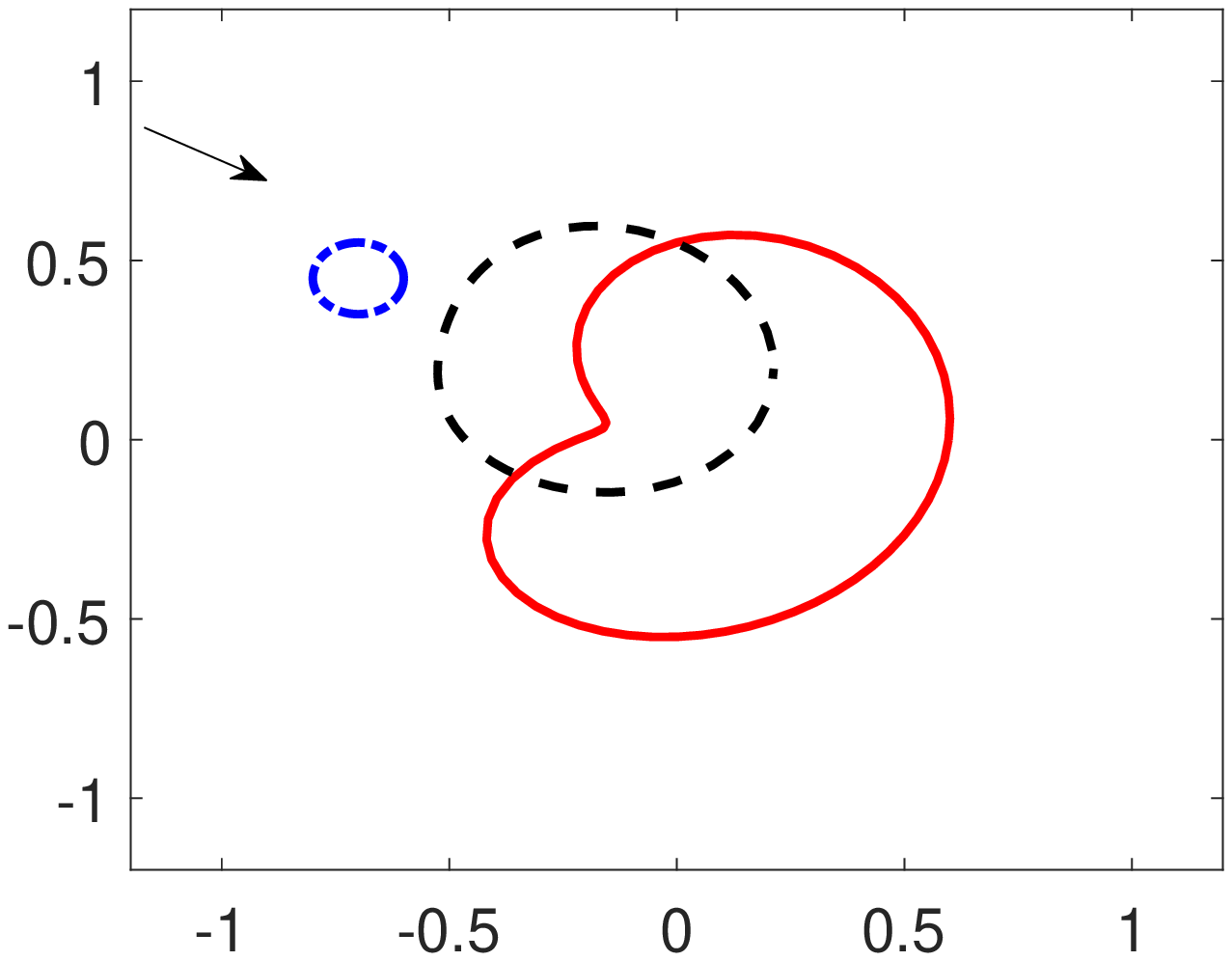}}
	\subfigure[10th iteration]{\includegraphics[width=0.32\textwidth]{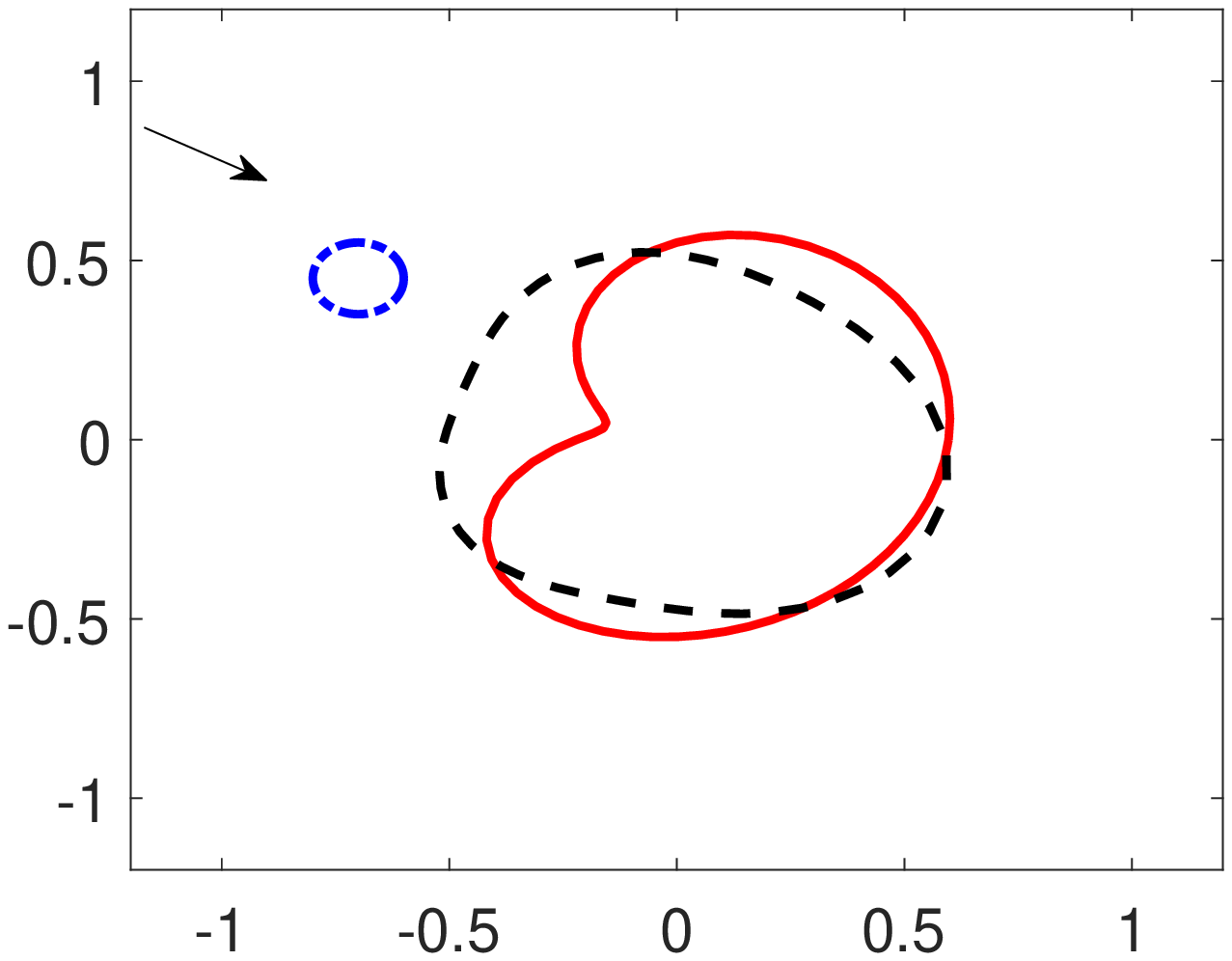}}
	\subfigure[15th iteration]{\includegraphics[width=0.32\textwidth]{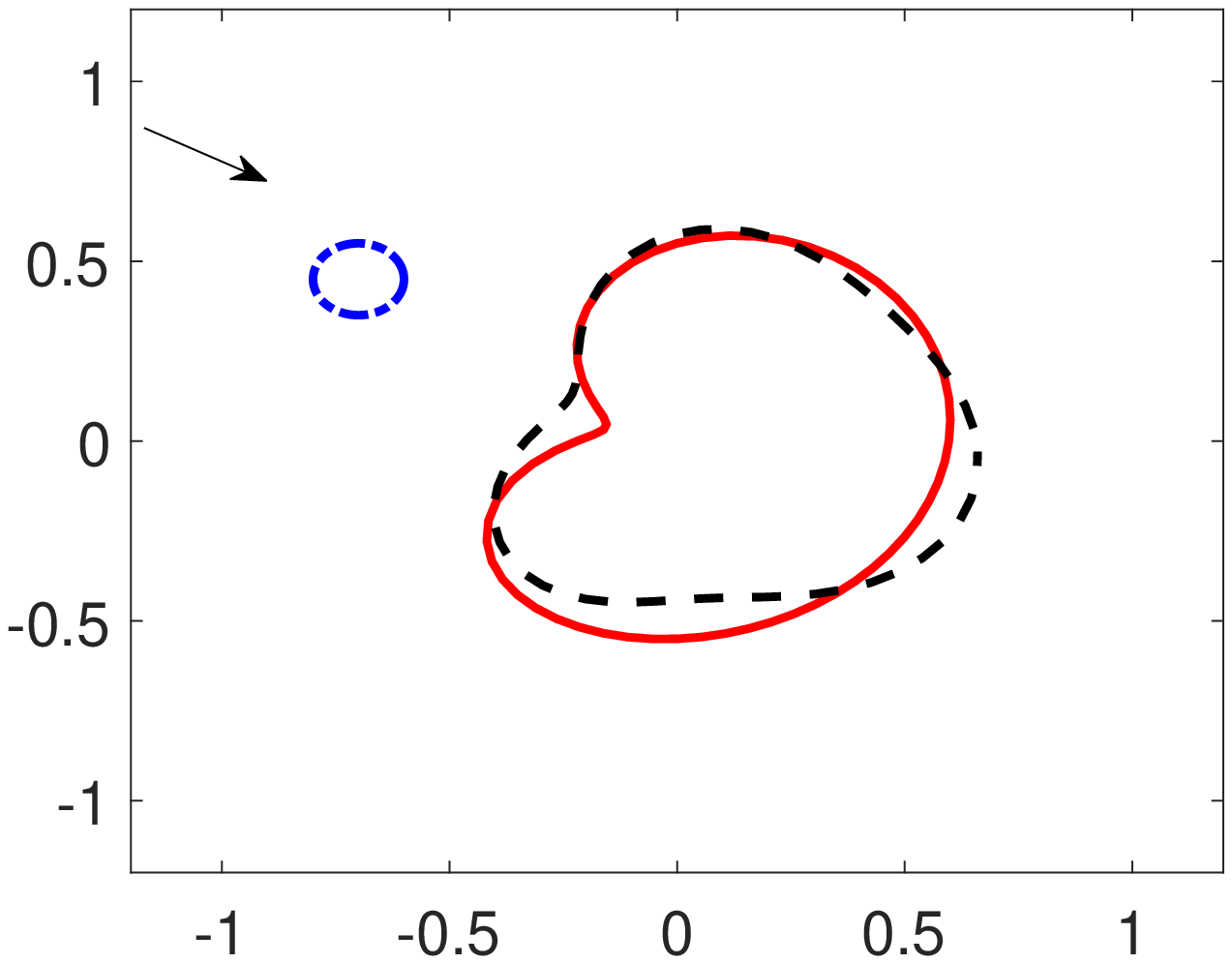}}
	\subfigure[20th iteration]{\includegraphics[width=0.32\textwidth]{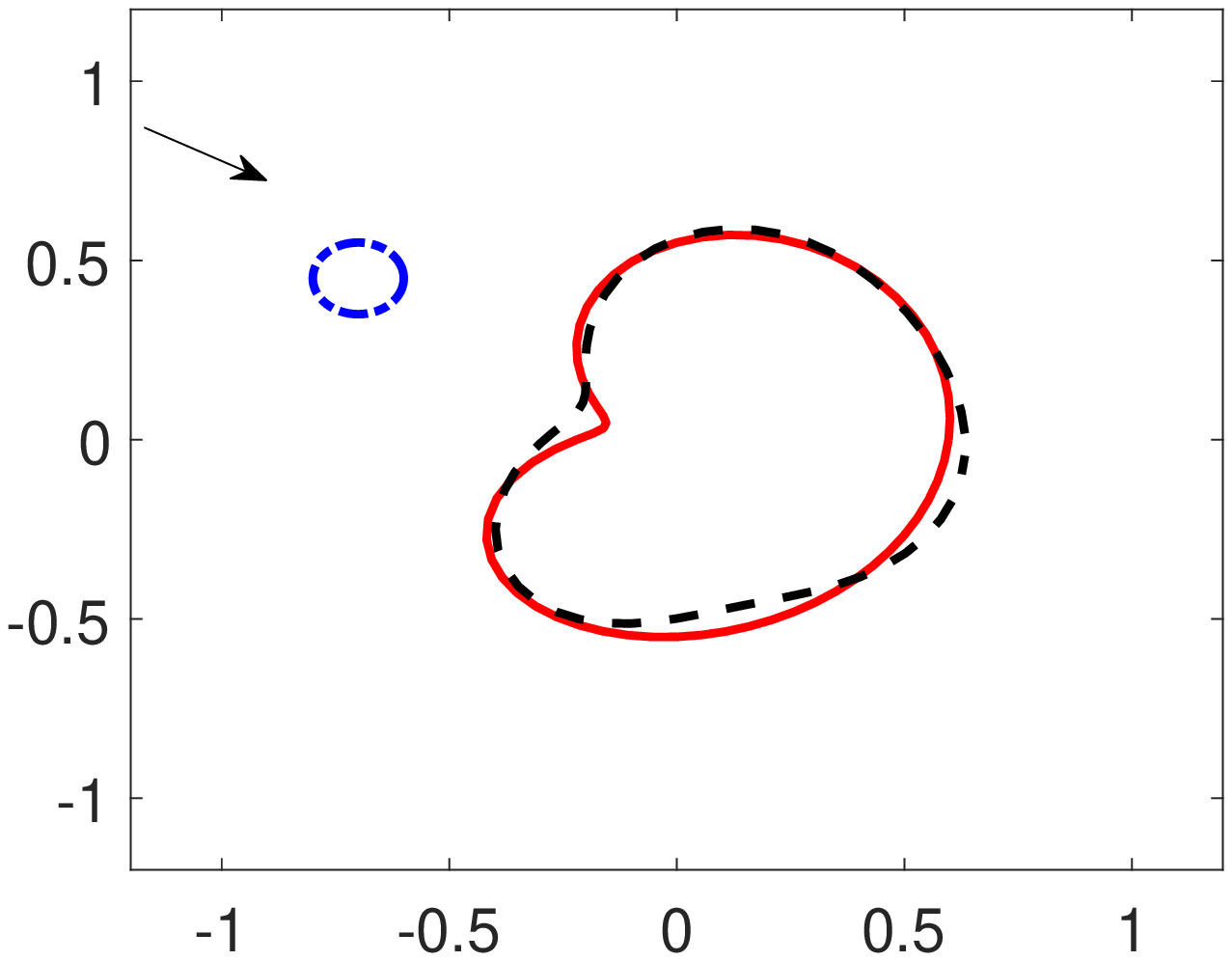}}
	\subfigure[24th iteration]{\includegraphics[width=0.32\textwidth]{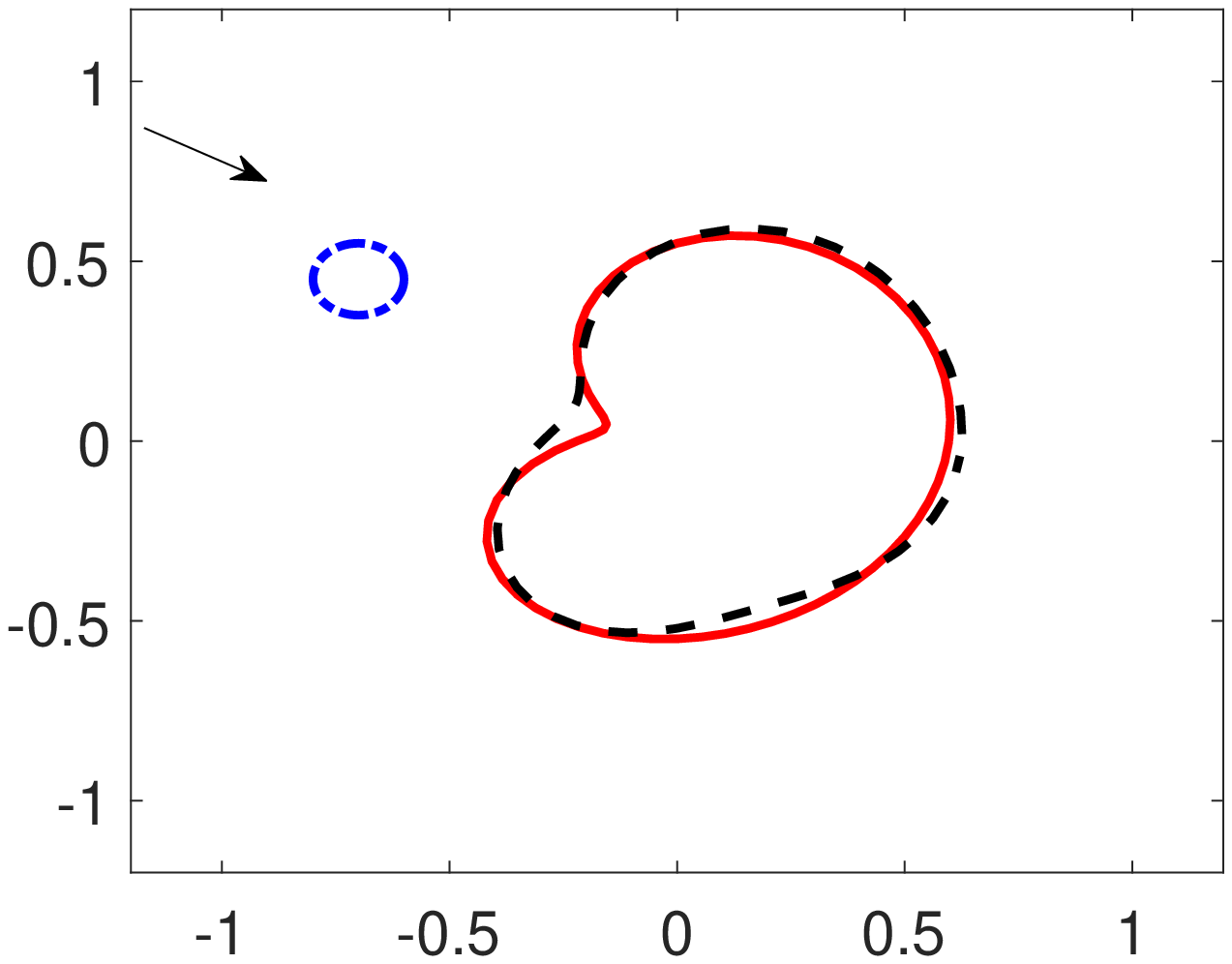}}
	\subfigure[relative error]{\includegraphics[width=0.32\textwidth]{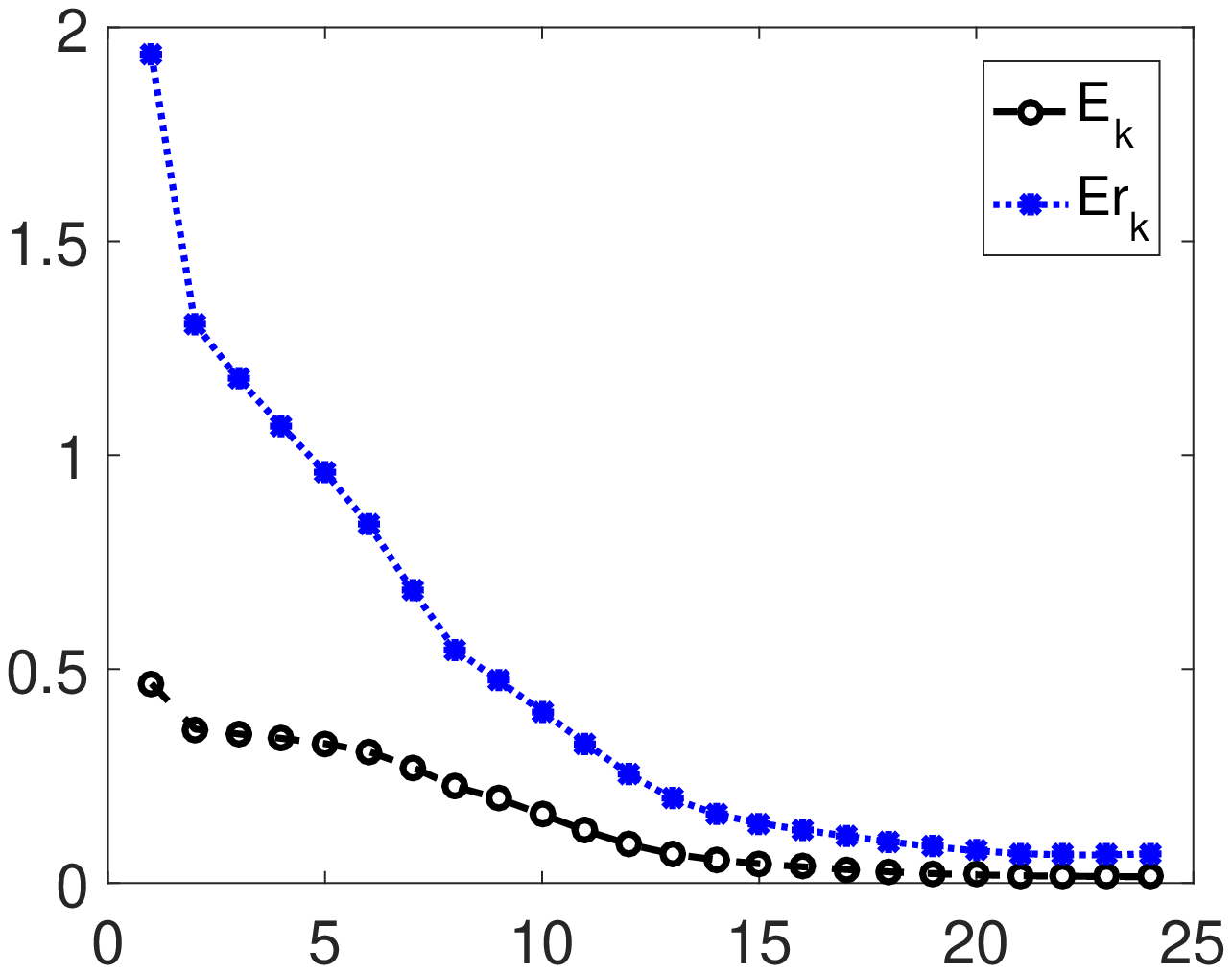}}
	\caption{Reconstructions of an apple-shaped domain with $1\%$ noise and $\epsilon=0.015$.}\label{4.1}
\end{figure}

\begin{figure}
	\centering 
	\subfigure[5th iteration]{\includegraphics[width=0.32\textwidth]{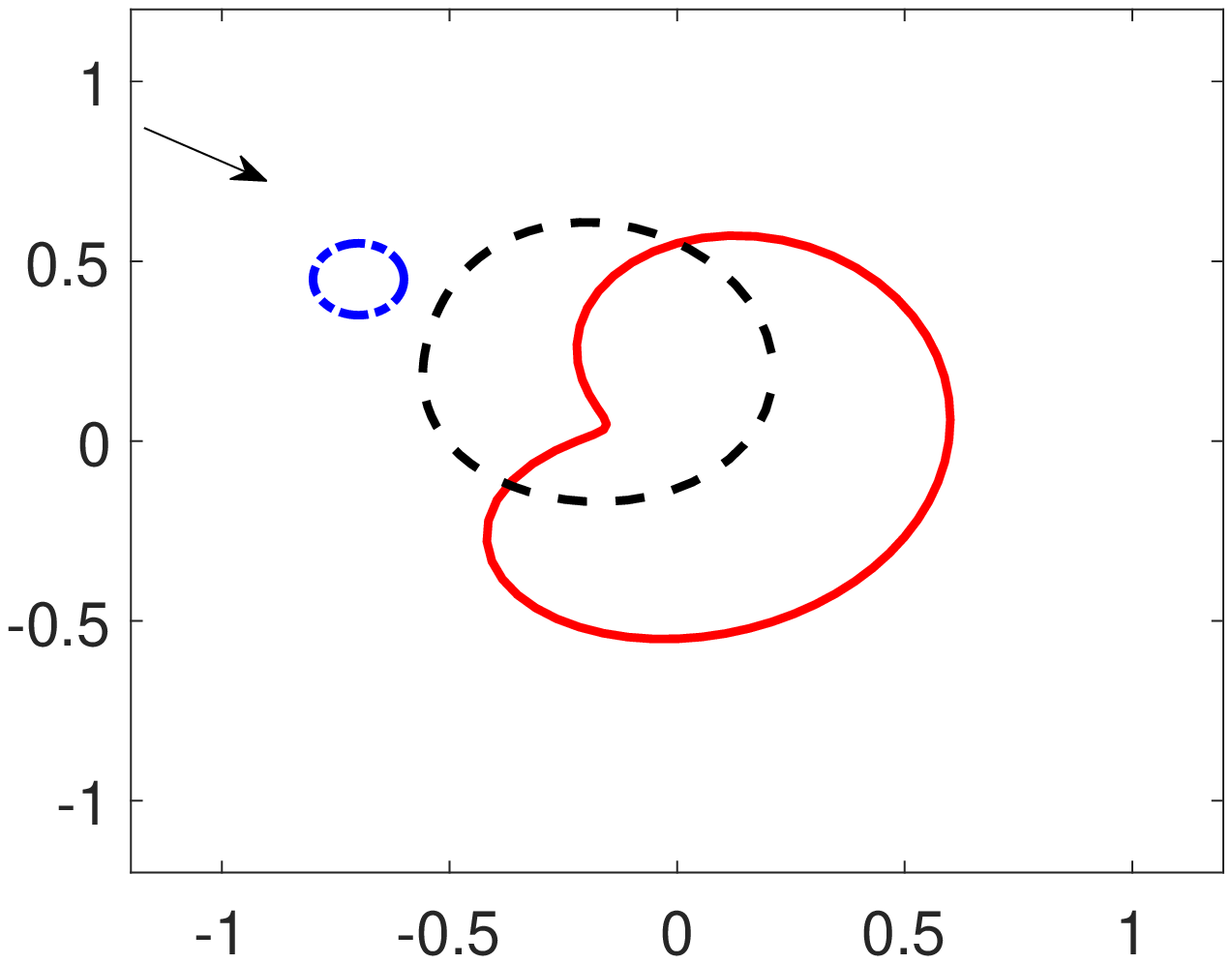}}
	\subfigure[10th iteration]{\includegraphics[width=0.32\textwidth]{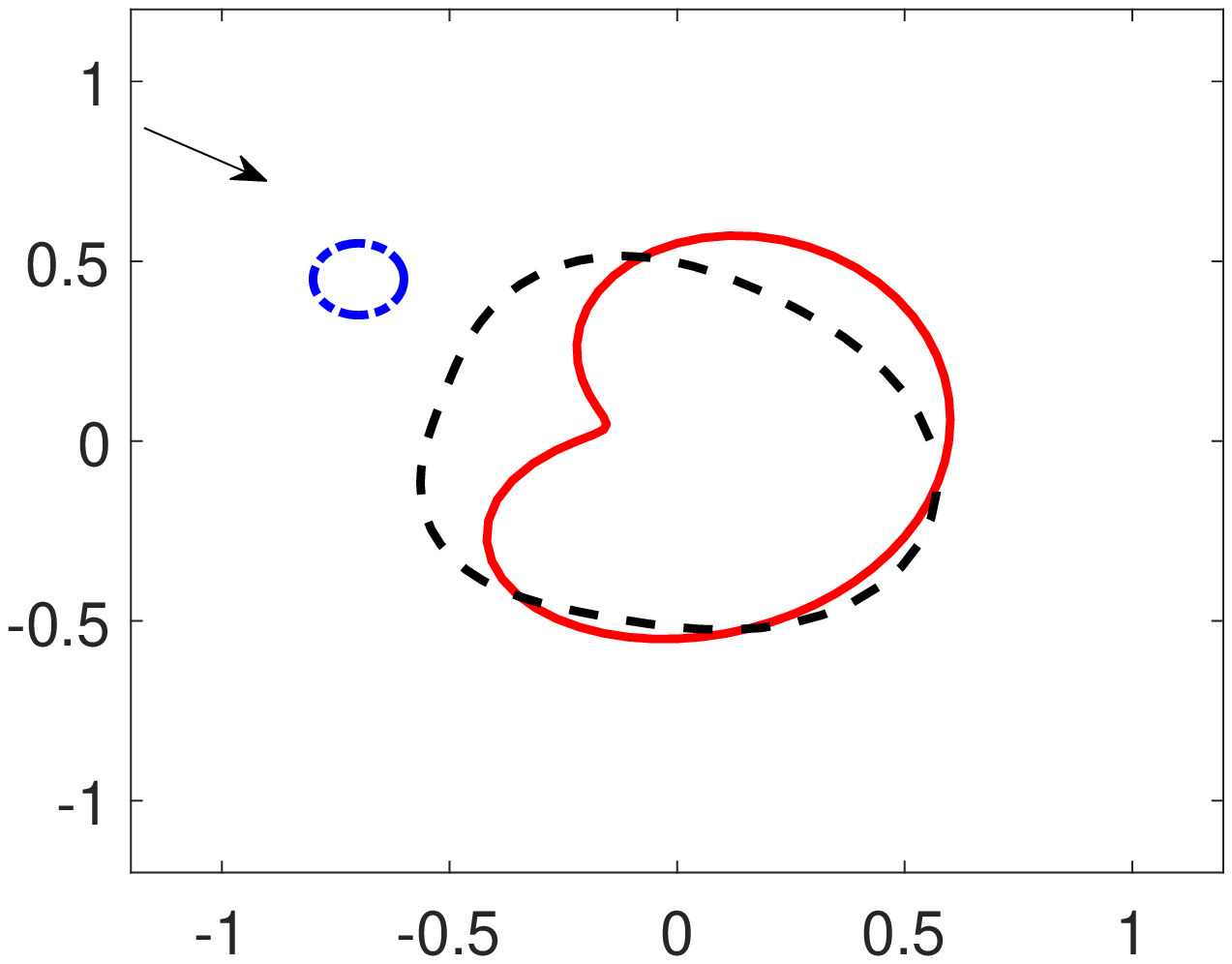}}
	\subfigure[15th iteration]{\includegraphics[width=0.32\textwidth]{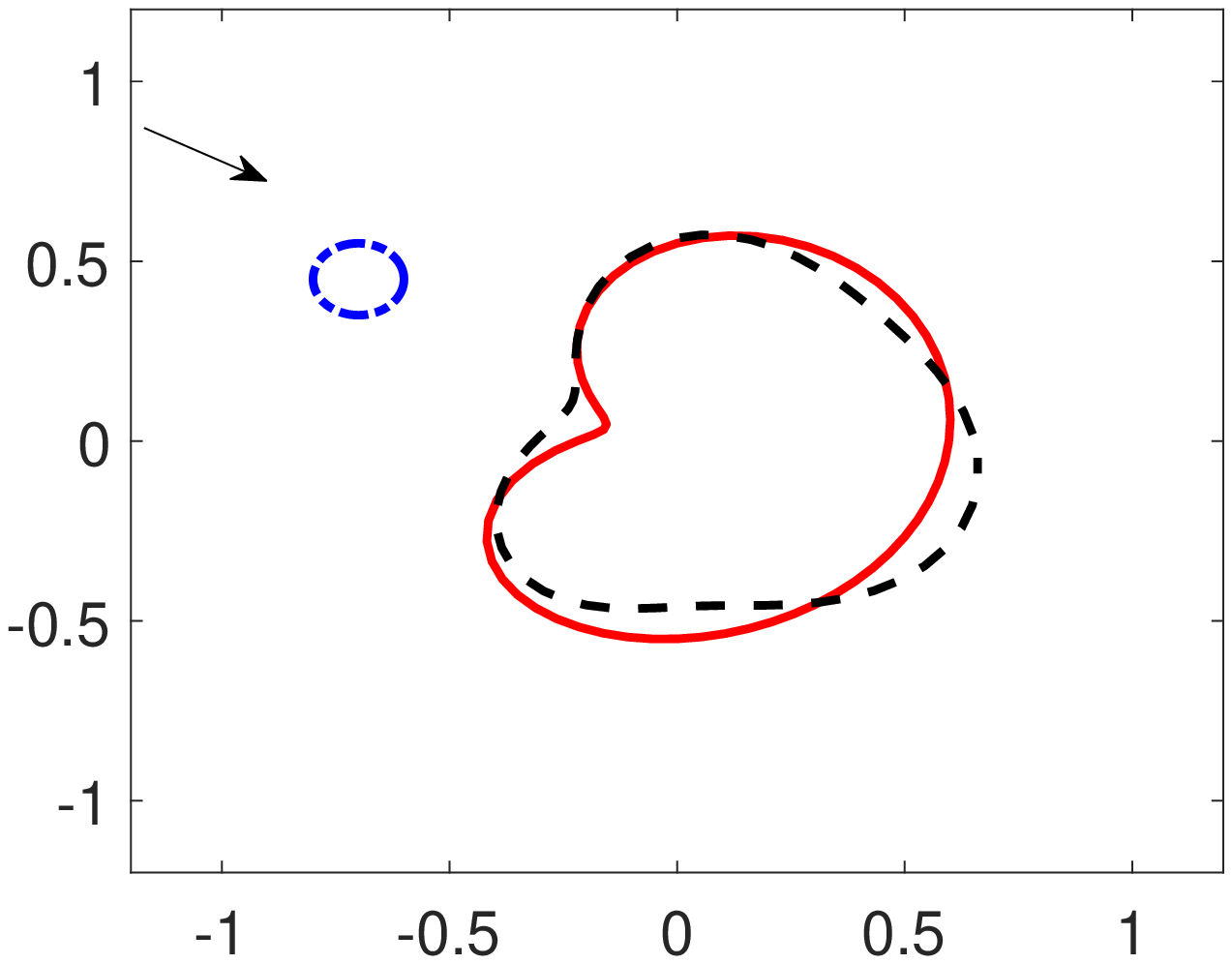}}
	\subfigure[18th iteration]{\includegraphics[width=0.32\textwidth]{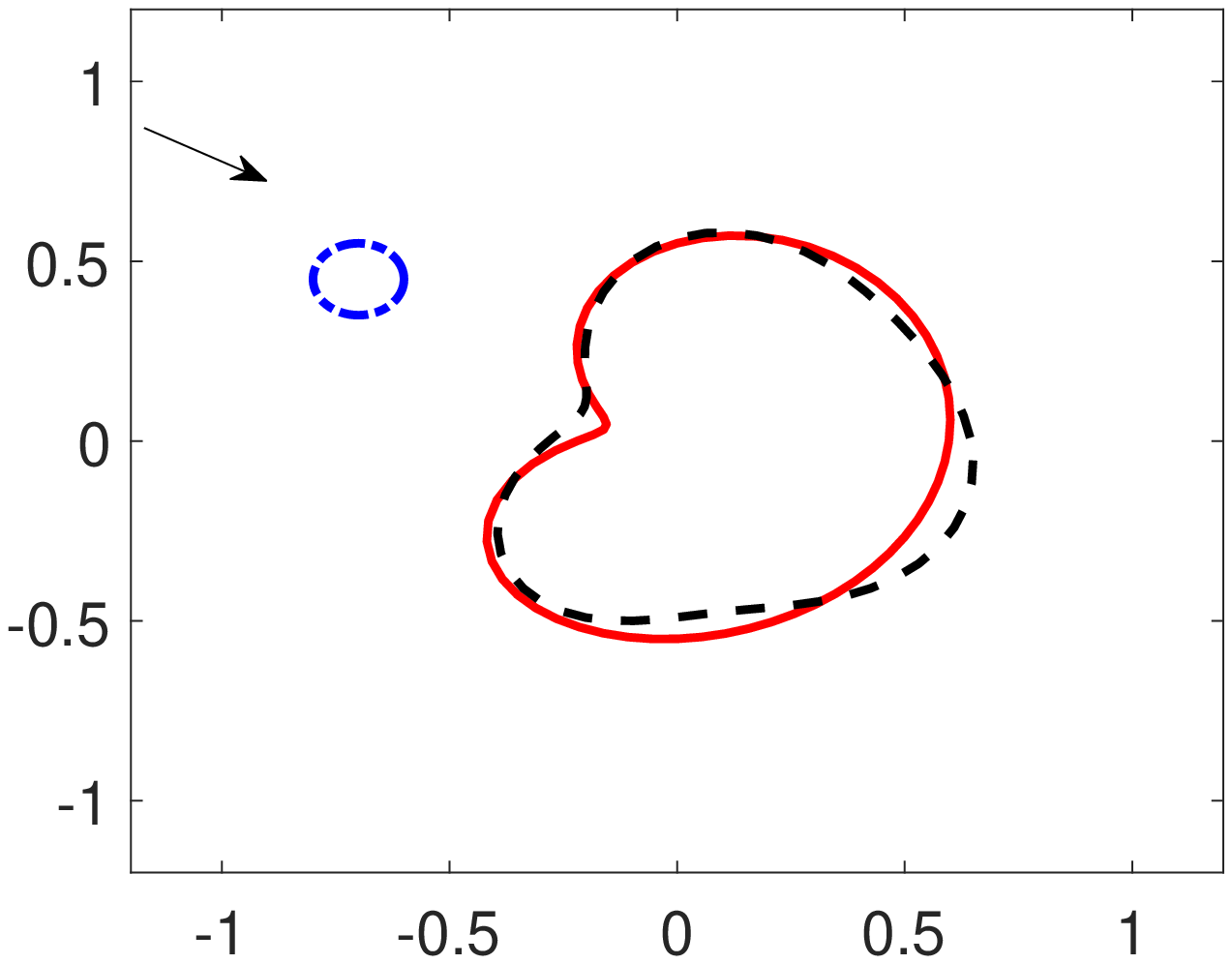}}
	\subfigure[relative error]{\includegraphics[width=0.32\textwidth]{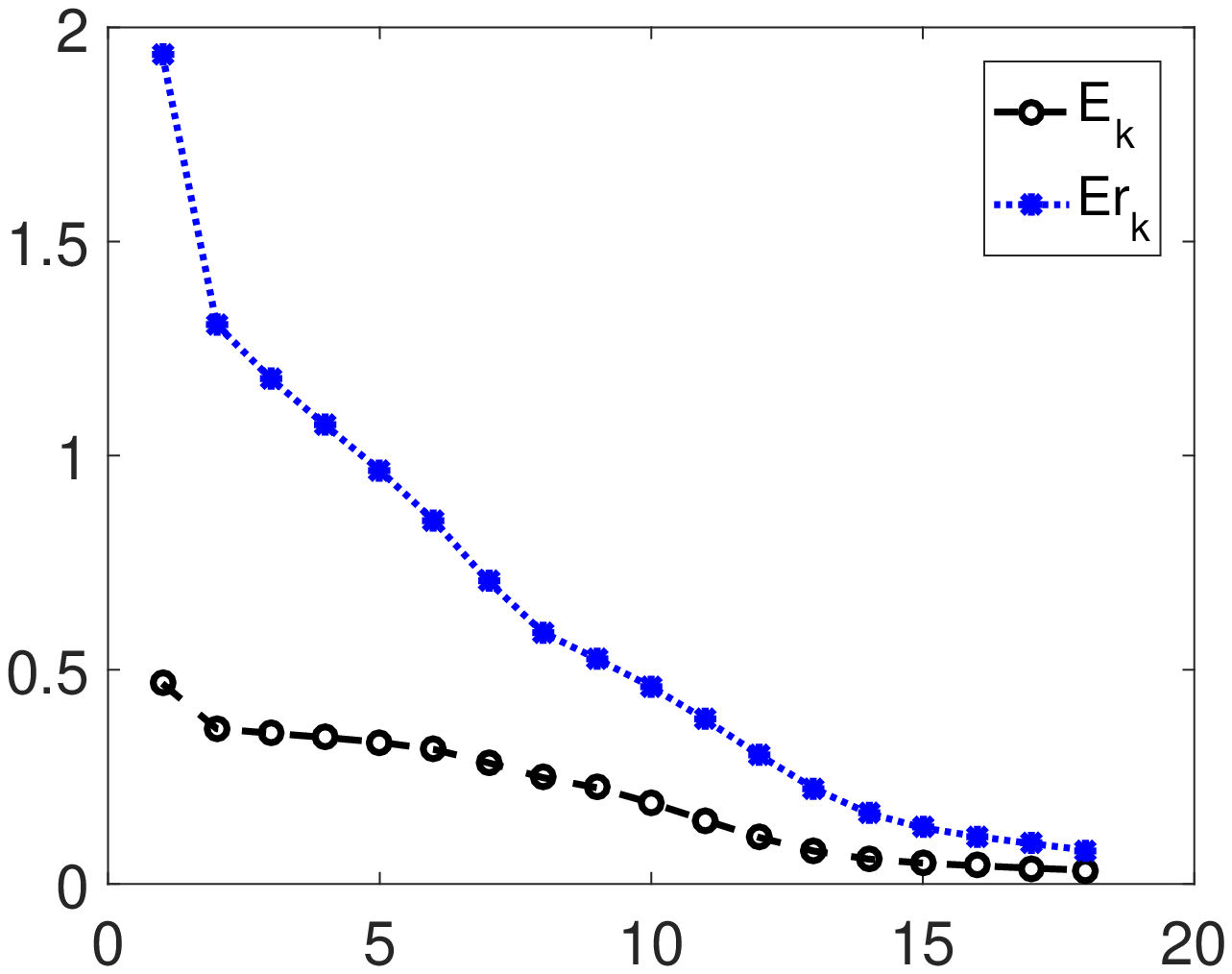}}
	\caption{Reconstructions of an apple-shaped domain with $5\%$ noise and $\epsilon=0.035$.}\label{4.2}
\end{figure}

\begin{figure}
	\centering 
	\subfigure[$(b_1, b_2)=(4, 0), R=0.4$]{\includegraphics[width=0.45\textwidth]{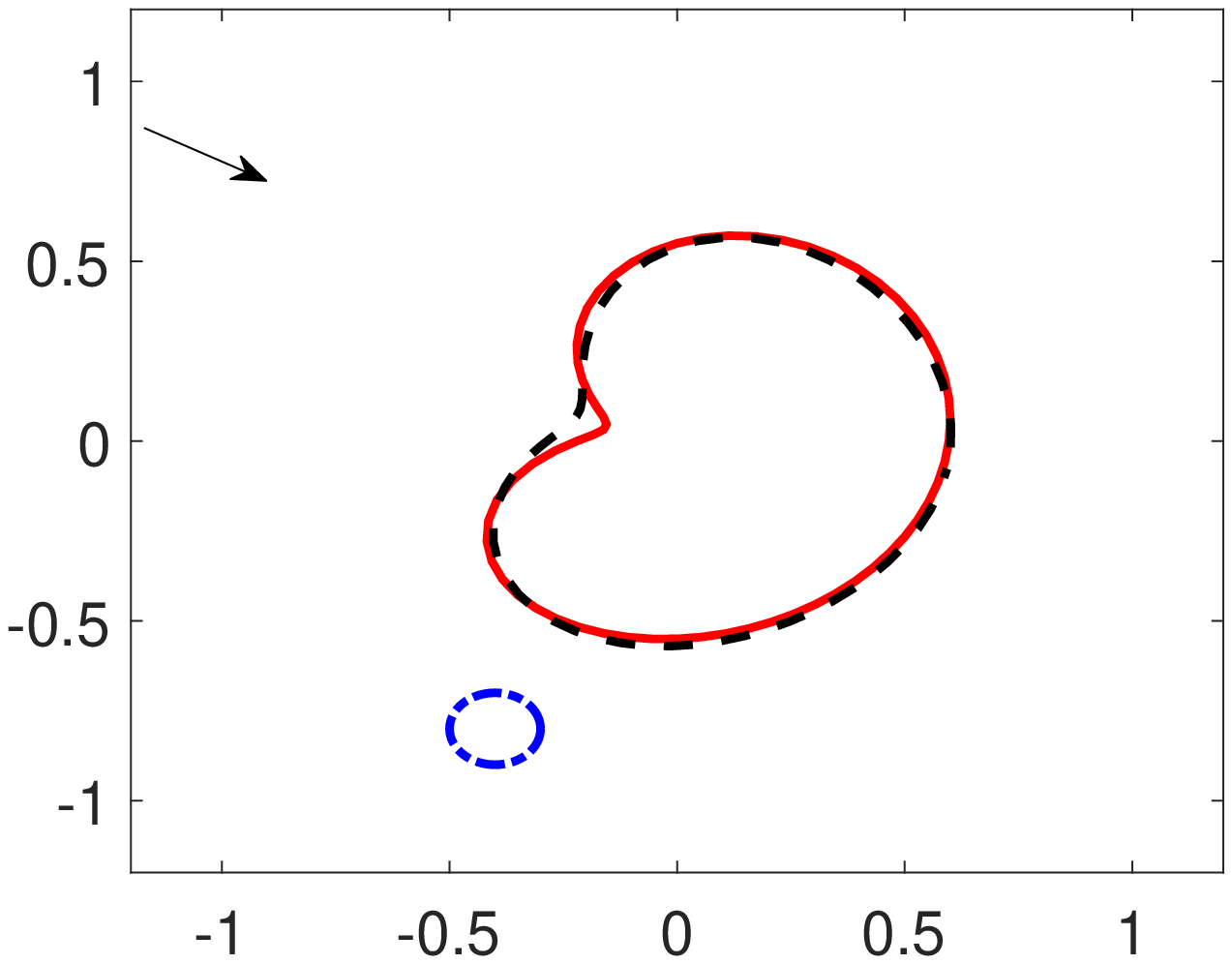}}
	\subfigure[$(b_1, b_2)=(4, 0), R=0.8$]{\includegraphics[width=0.45\textwidth]{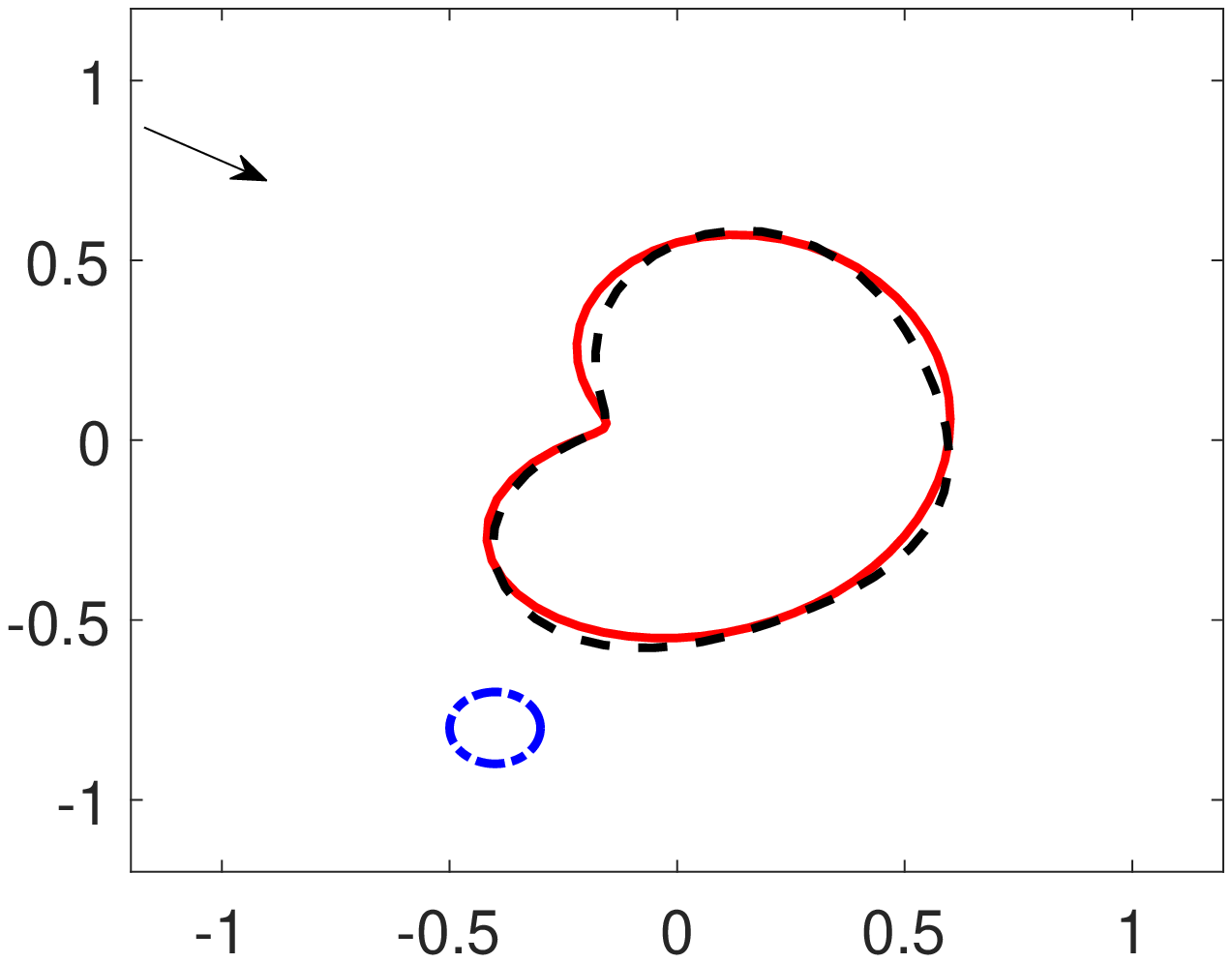}} 
	\subfigure[$(b_1, b_2)=(6, 0), R=0.4$]{\includegraphics[width=0.45\textwidth]{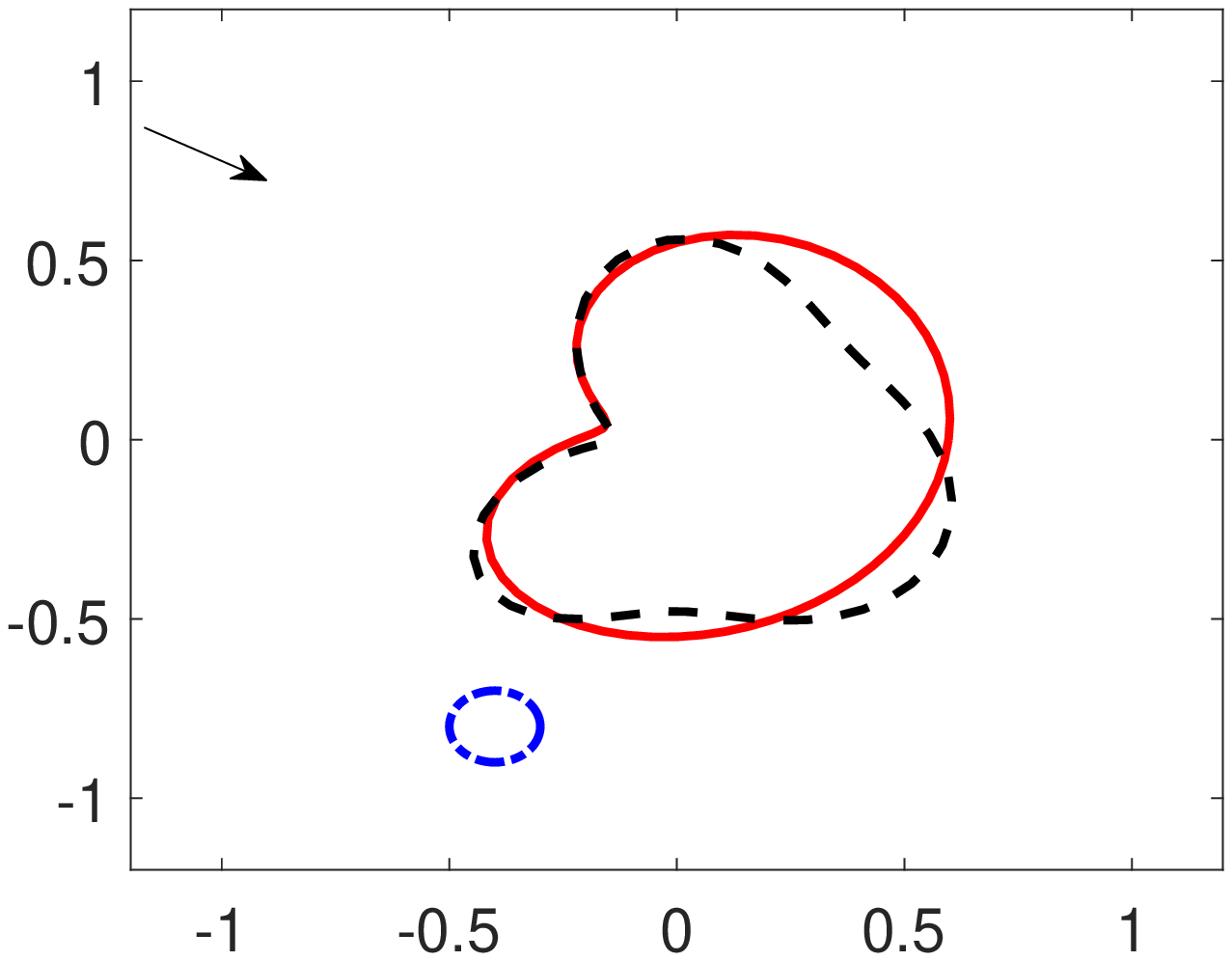}}
	\subfigure[$(b_1,b_2)=(6, 0), R=0.8$]{\includegraphics[width=0.45\textwidth]{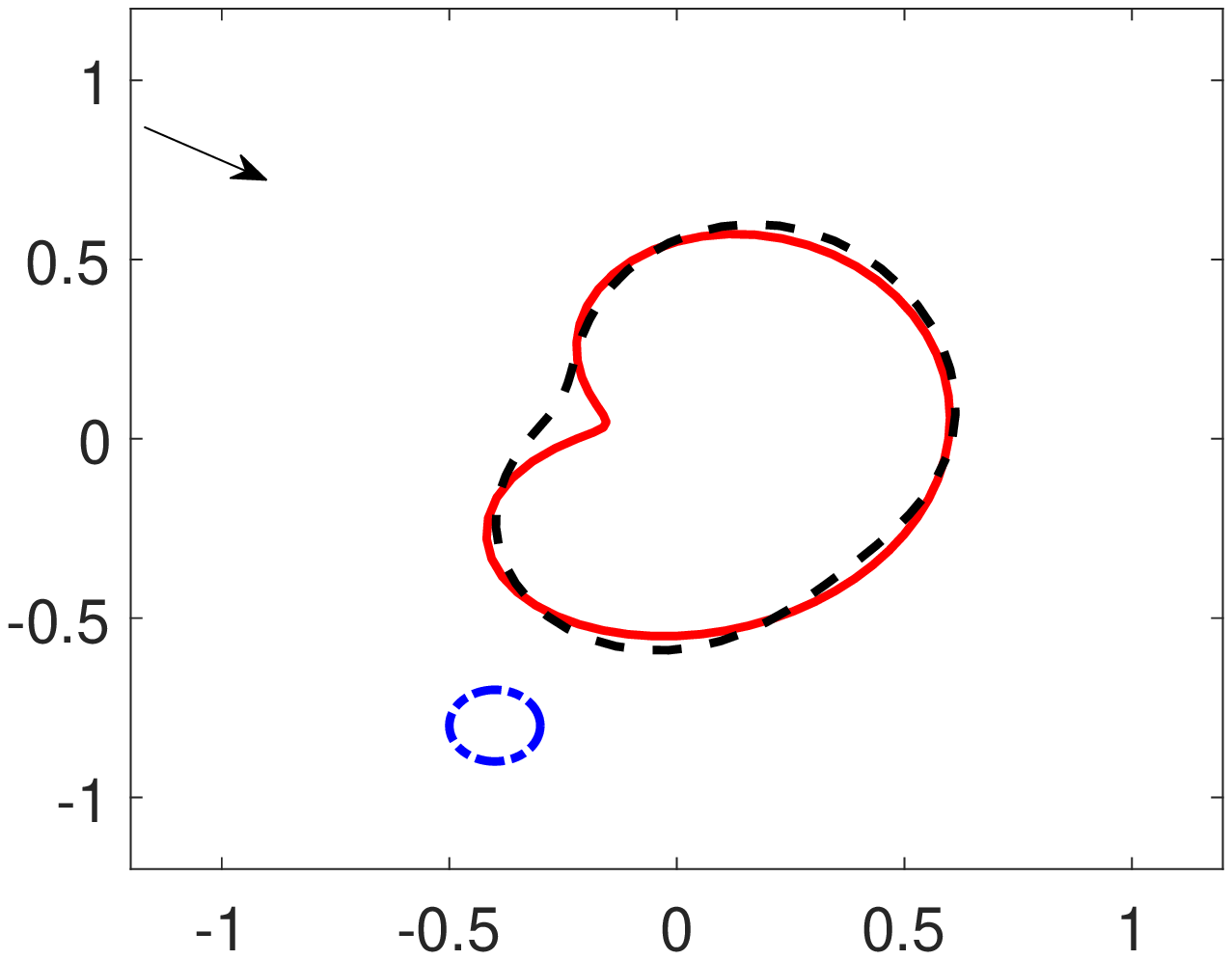}}
	\caption{Reconstructions of an apple-shaped domain with different reference balls and $1\%$ noise. Here we are using  the initial guess $(c_1^{(0)},c_2^{(0)})=(-0.4, -0.8), r^{(0)}=0.1$, and (a) $\epsilon=0.011$, (b) $\epsilon=0.02$, (c) $\epsilon=0.03$, (d) $\epsilon=0.011$.}\label{4.3}
\end{figure}

\begin{figure}
	\centering 
	\subfigure[$d=(\cos(-\pi/6),\sin(-\pi/6))$]{\includegraphics[width=0.45\textwidth]{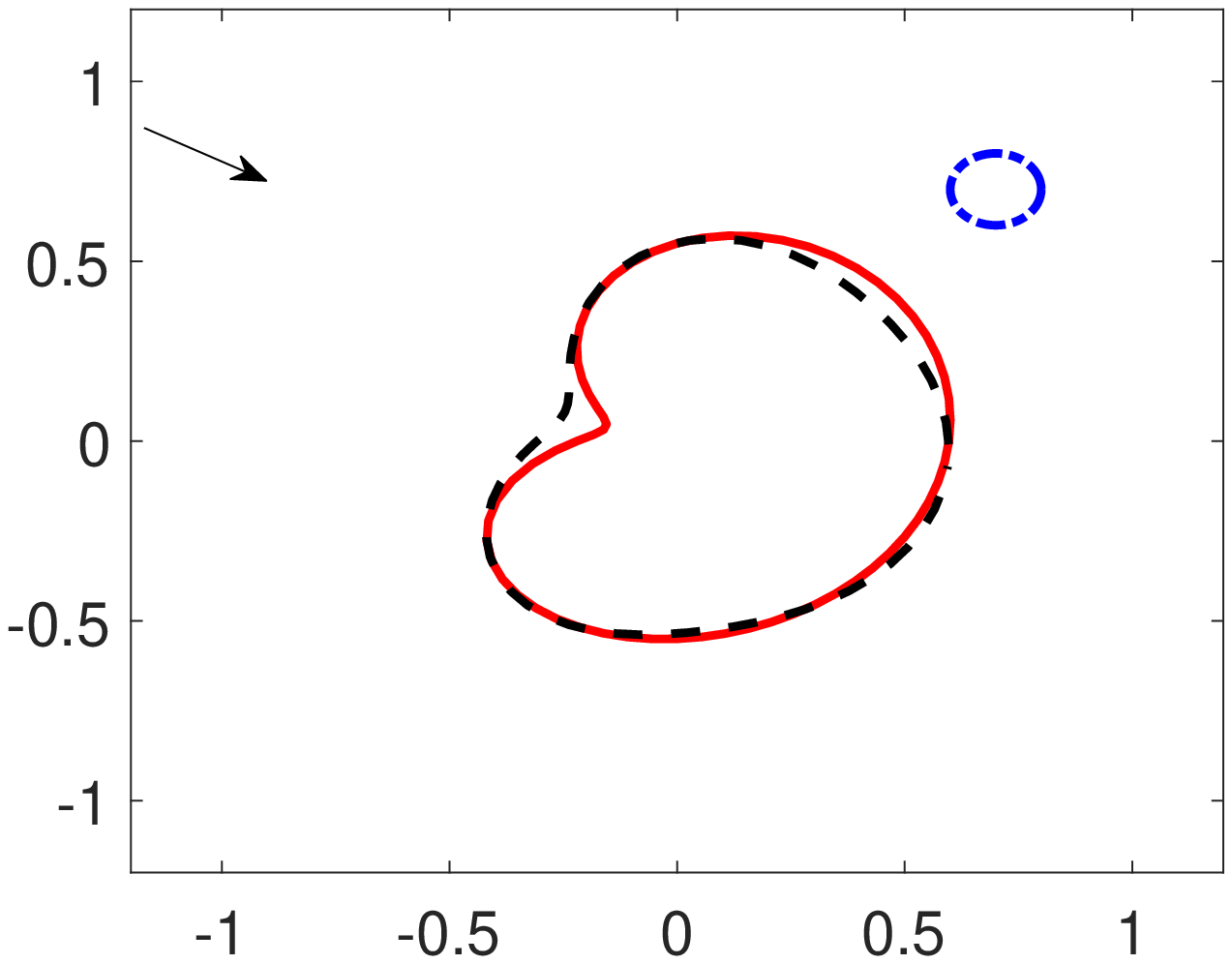}}
	\subfigure[$d=(\cos(4\pi/3),\sin(4\pi/3))$]{\includegraphics[width=0.45\textwidth]{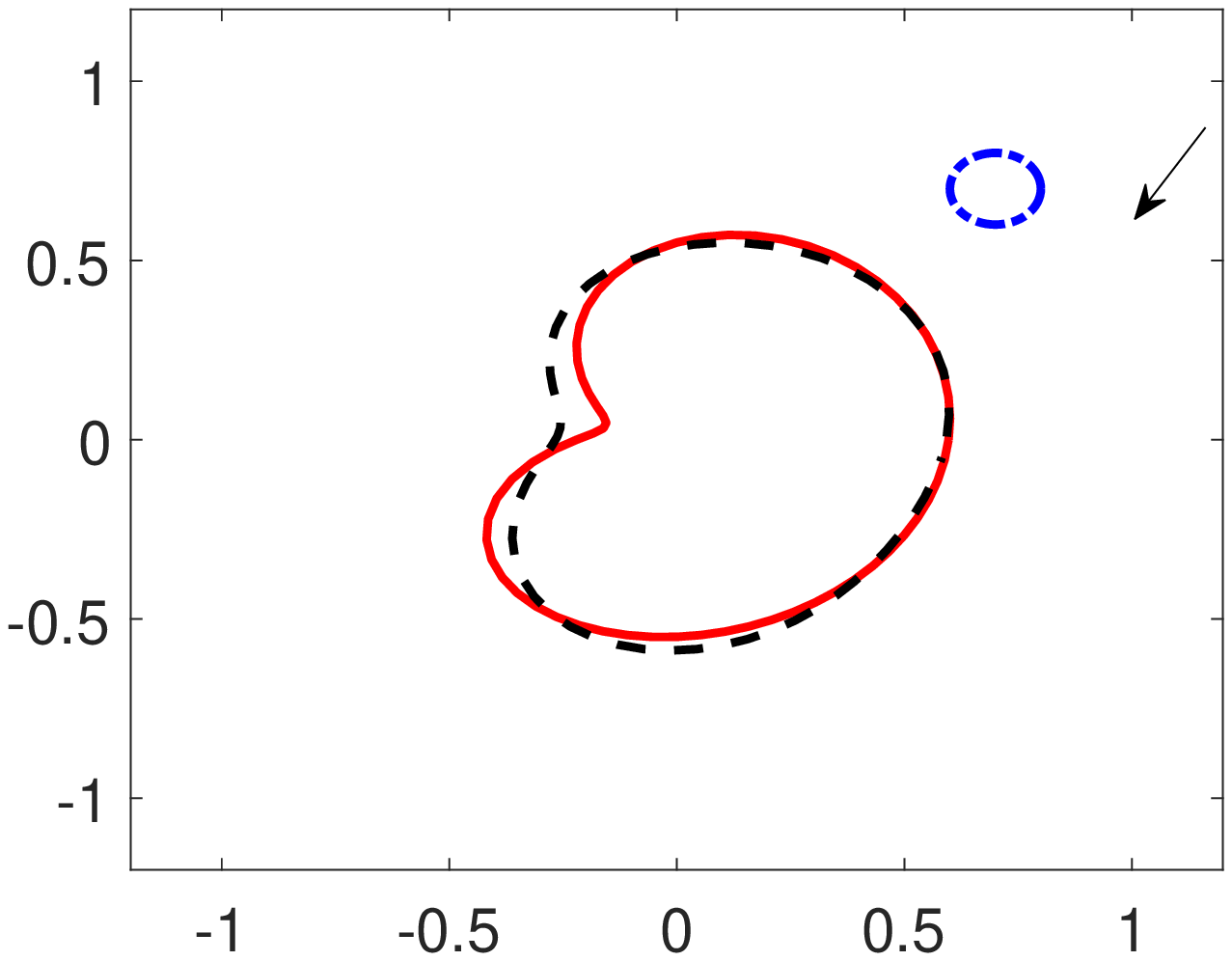}}
	\caption{Reconstructions of an apple-shaped domain with different incoming wave directions, $1\%$ noise is added. Here, the initial guess $(c_1^{(0)},c_2^{(0)})=(0.7, 0.7), r^{(0)}=0.1$, the reference ball $(b_1, b_2)=(4, 0), R=0.4$, and $\epsilon=0.015$. }\label{4.4}
\end{figure}

\begin{figure}
	\centering \subfigure[$(b_1, b_2)=(4, 0), R=0.2$]{\includegraphics[width=0.45\textwidth]{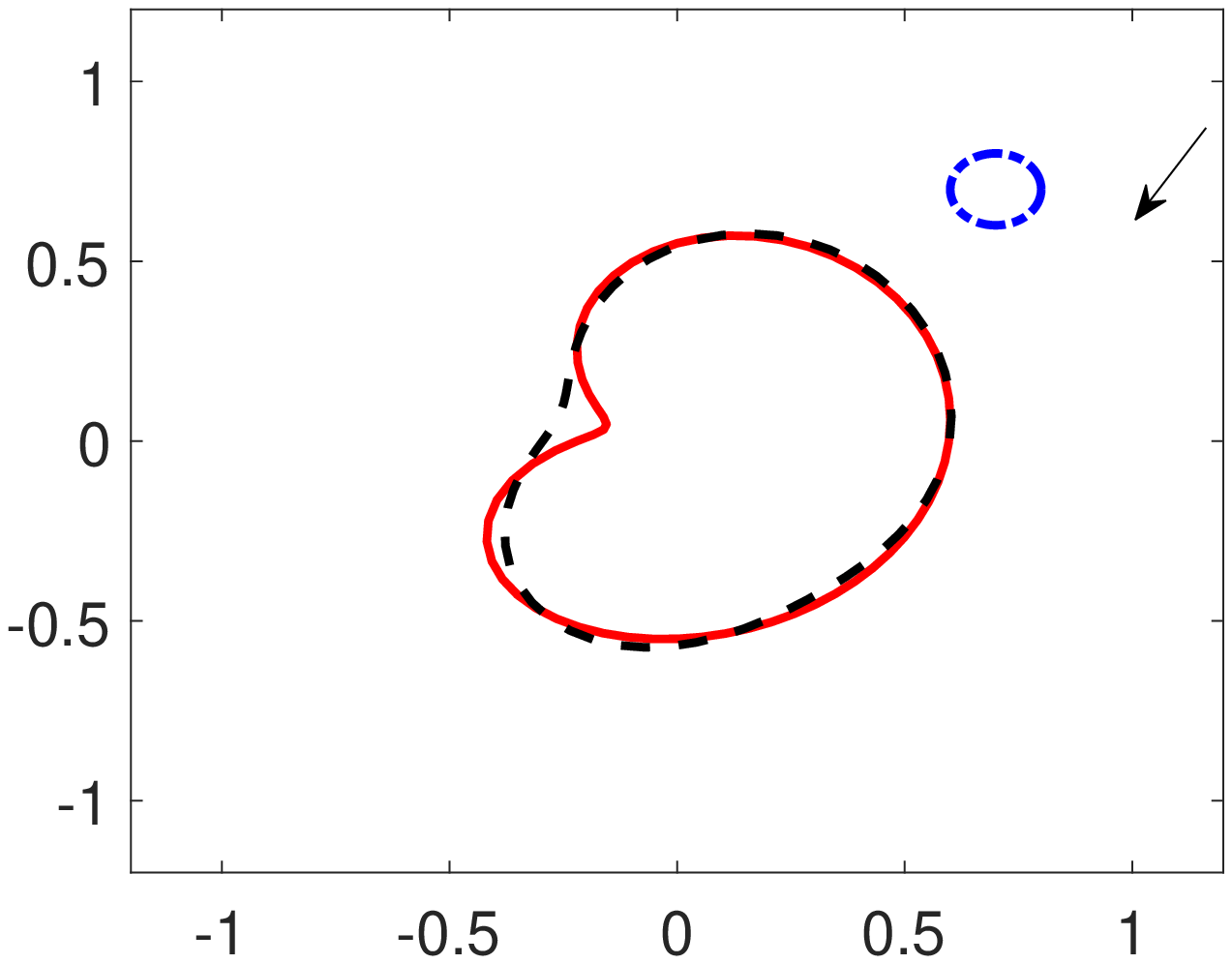}} 
	\subfigure[$(b_1, b_2)=(6, 0), R=0.2$]{\includegraphics[width=0.45\textwidth]{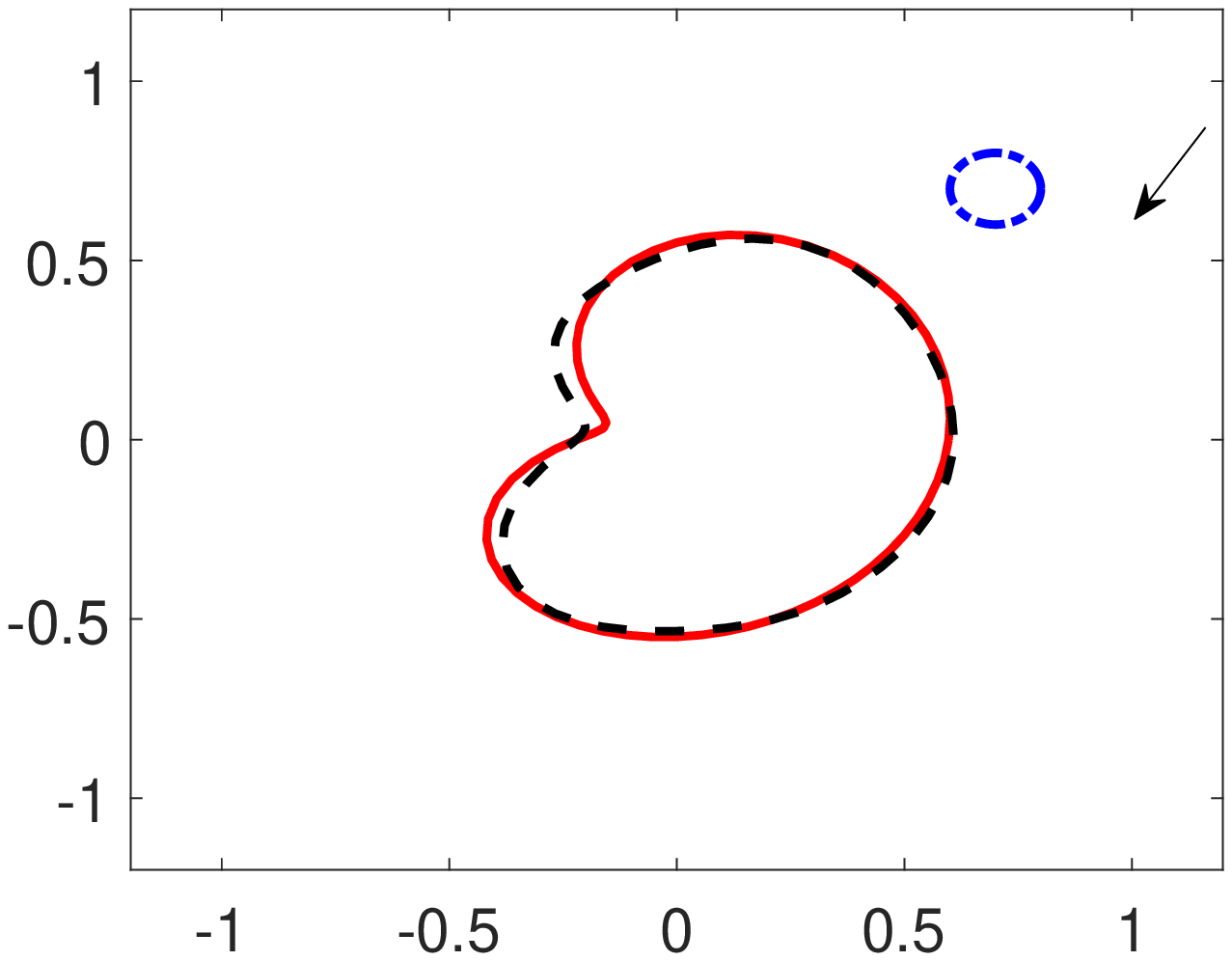}}
	\caption{Reconstructions of an apple-shaped domain with different locations of the reference ball, $1\%$ noise is added. Here, the initial guess $(c_1^{(0)},c_2^{(0)})=(0.7, 0.7), r^{(0)}=0.1$ and $\epsilon=0.015$. }\label{4.5}
\end{figure}

\begin{figure}
	\centering 
	\subfigure[5th iteration]{\includegraphics[width=0.32\textwidth]{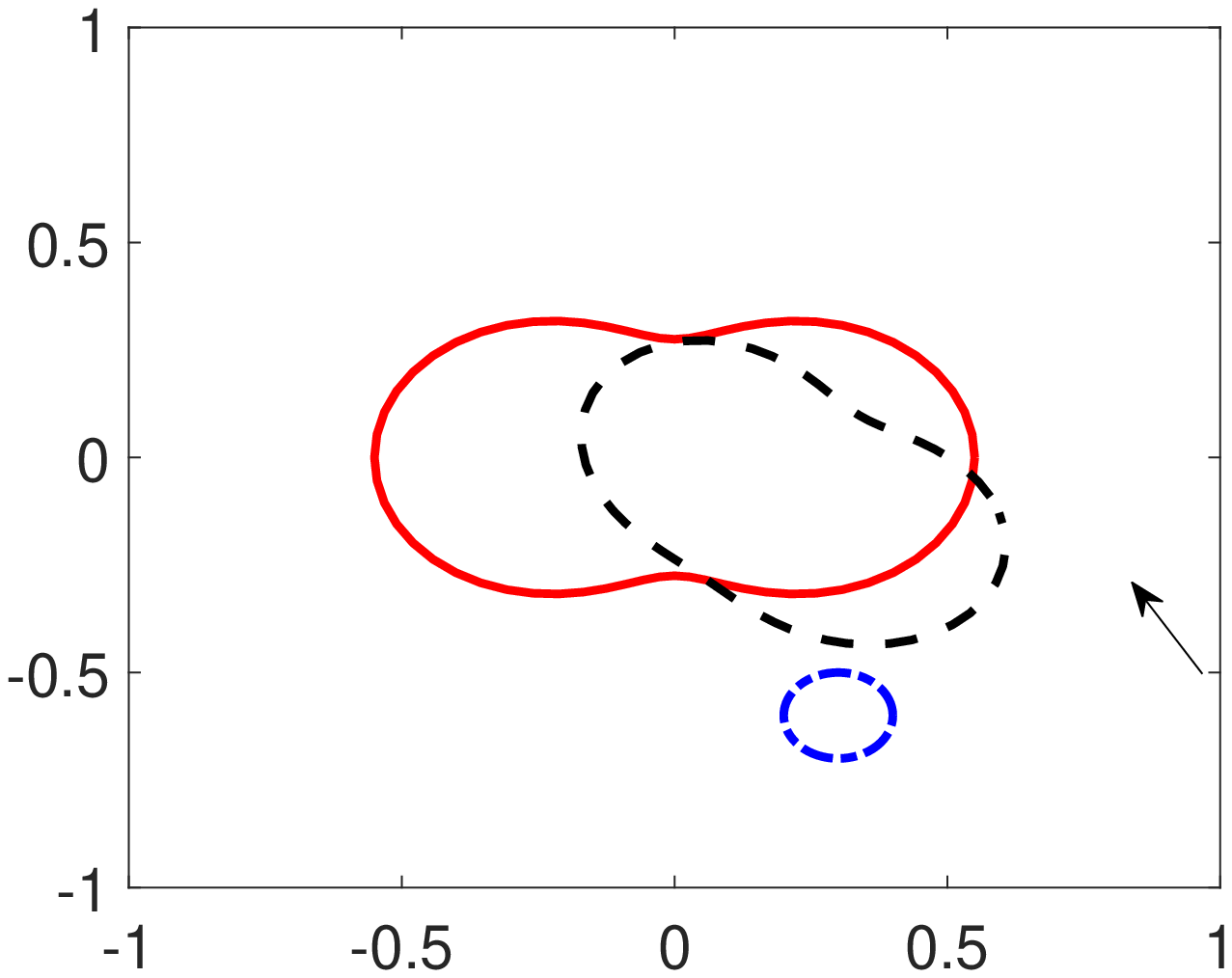}}
	\subfigure[10th iteration]{\includegraphics[width=0.32\textwidth]{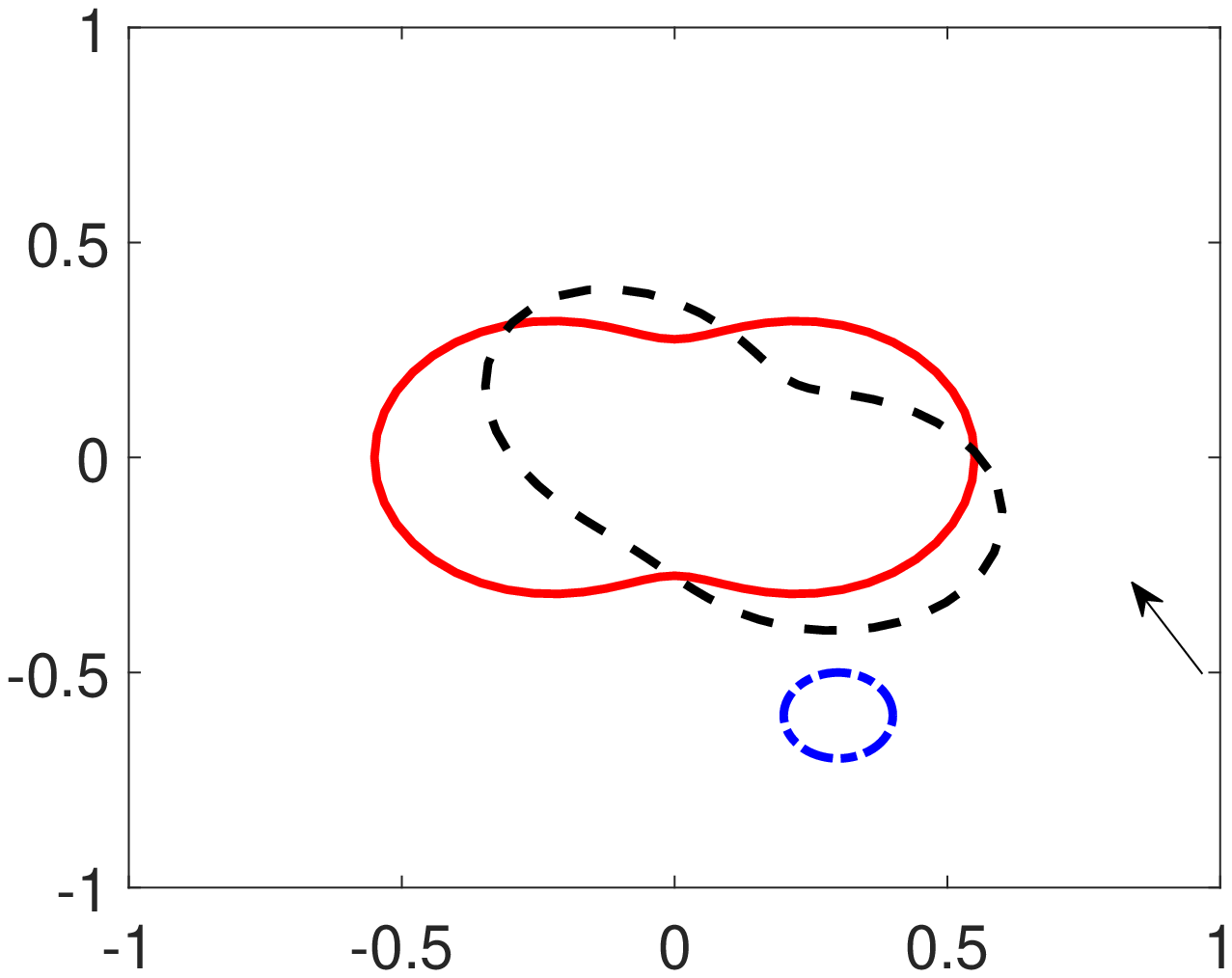}}
	\subfigure[15th iteration]{\includegraphics[width=0.32\textwidth]{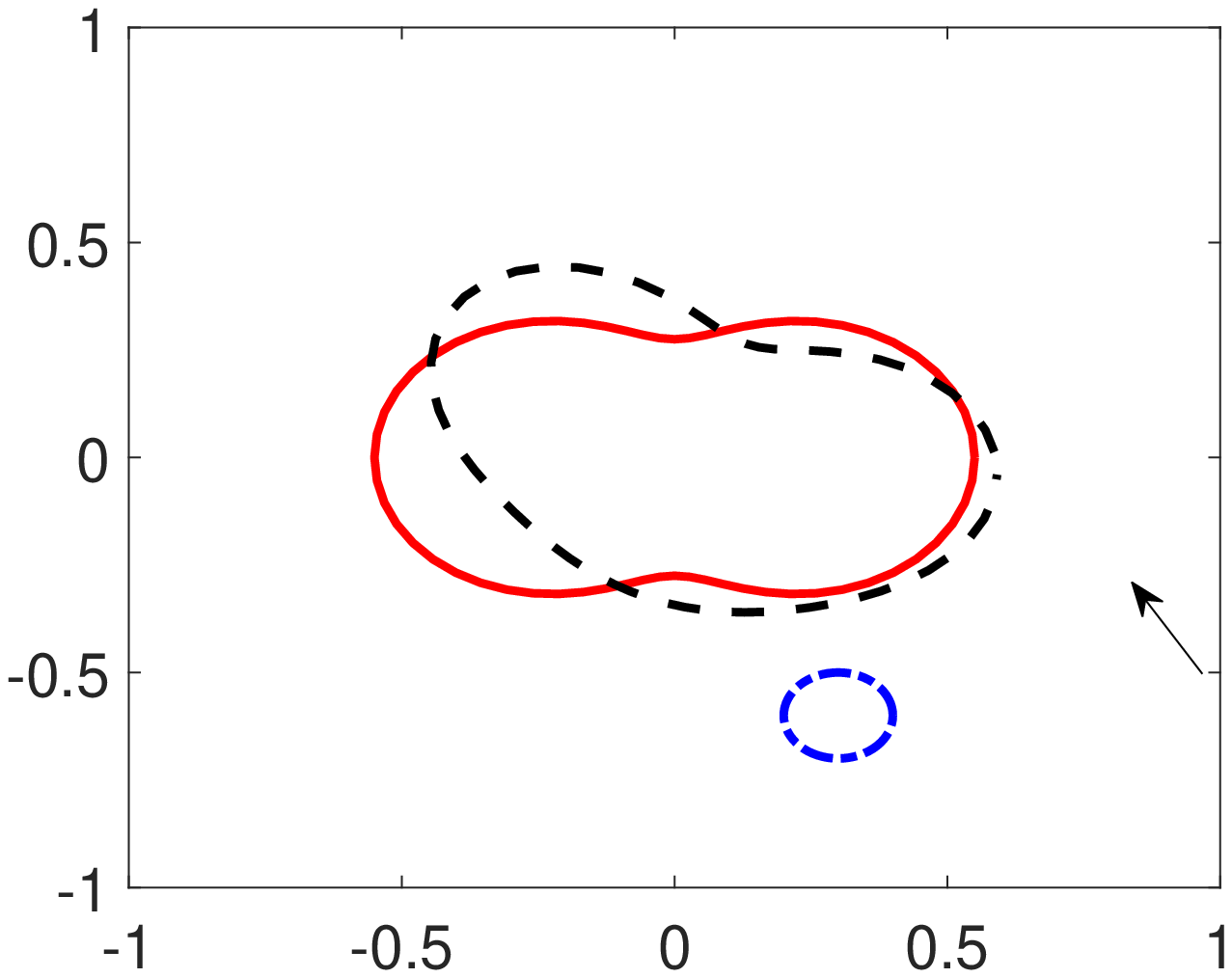}}
	\subfigure[20th iteration]{\includegraphics[width=0.32\textwidth]{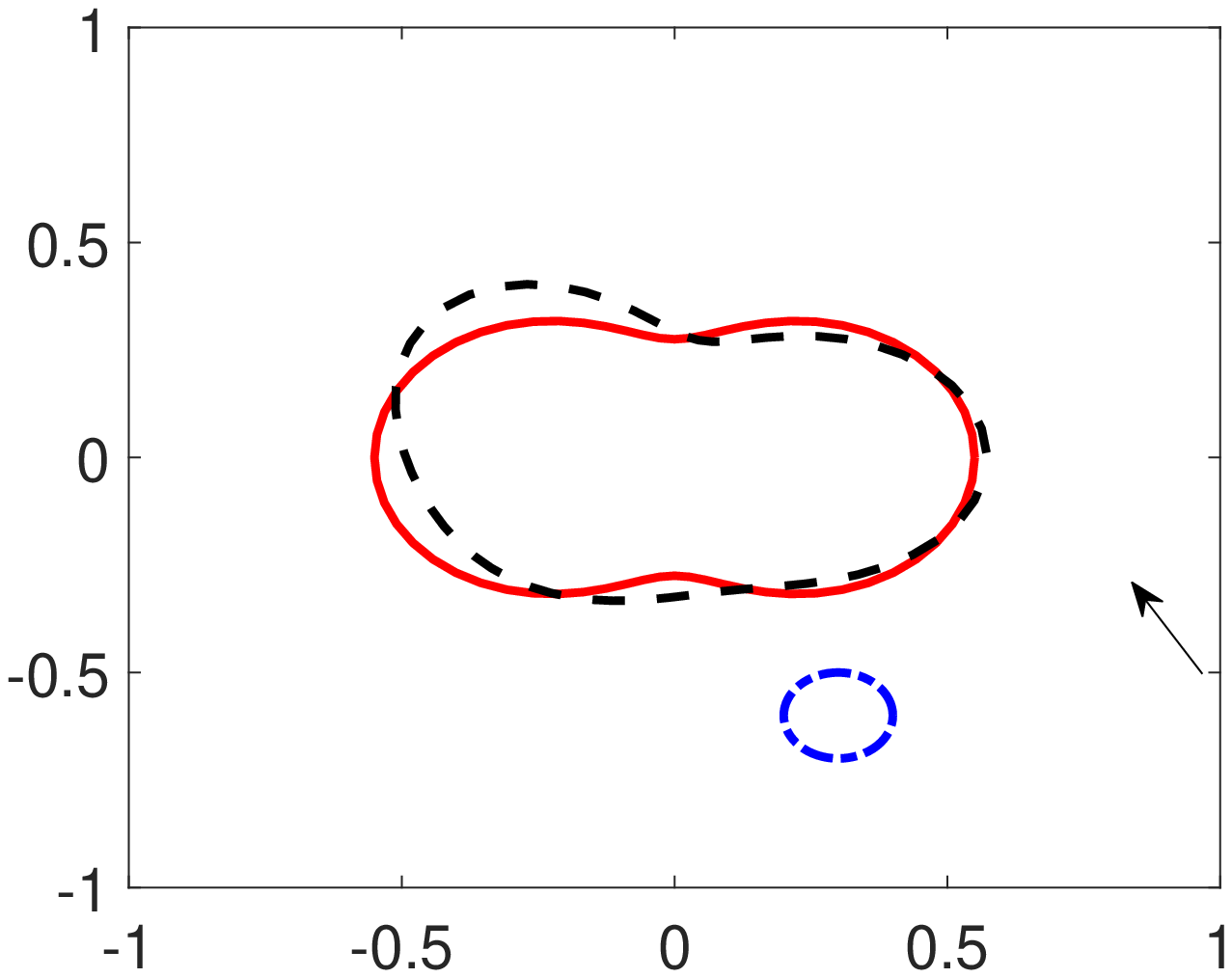}}
	\subfigure[23th iteration]{\includegraphics[width=0.32\textwidth]{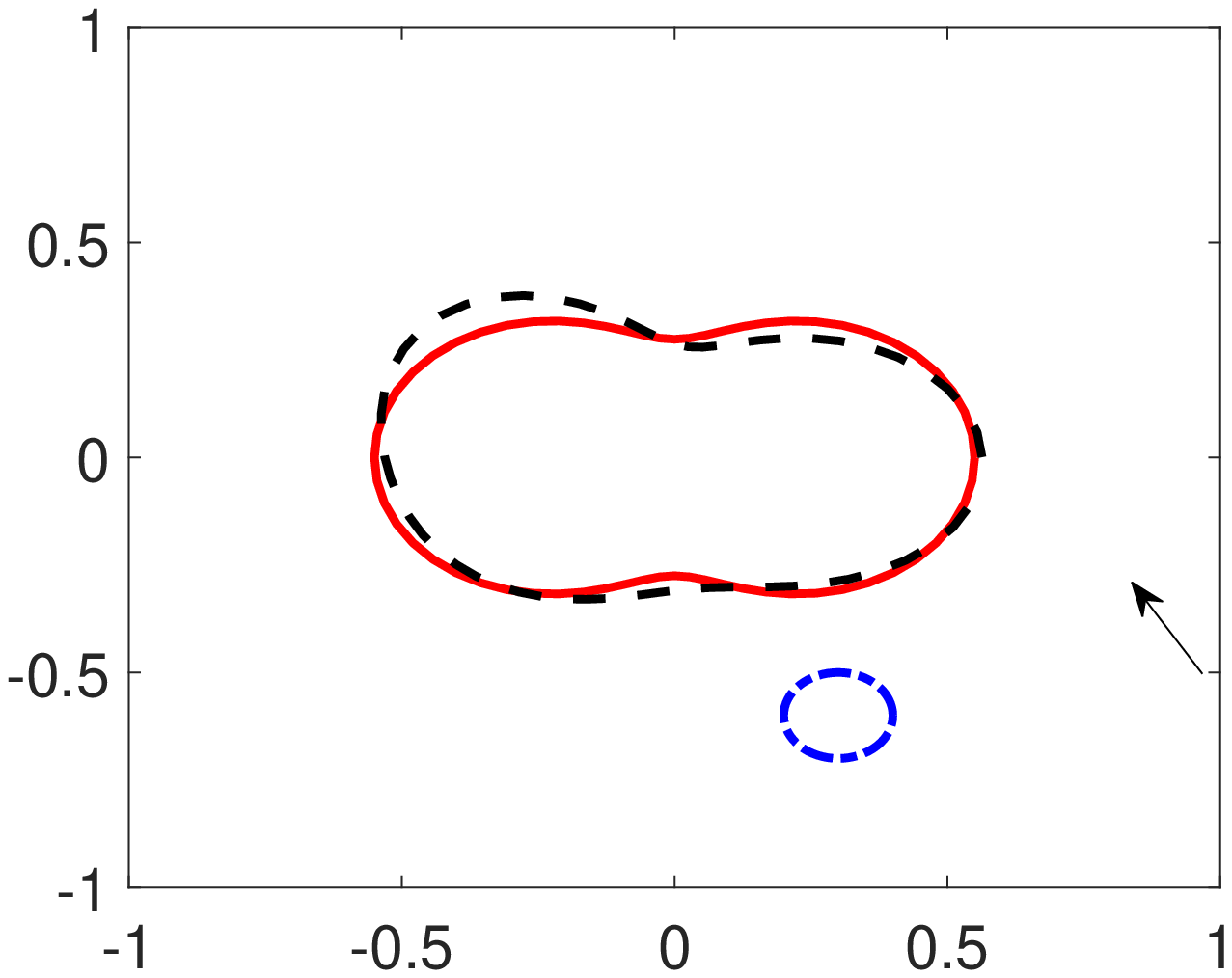}}
	\subfigure[relative error]{\includegraphics[width=0.32\textwidth]{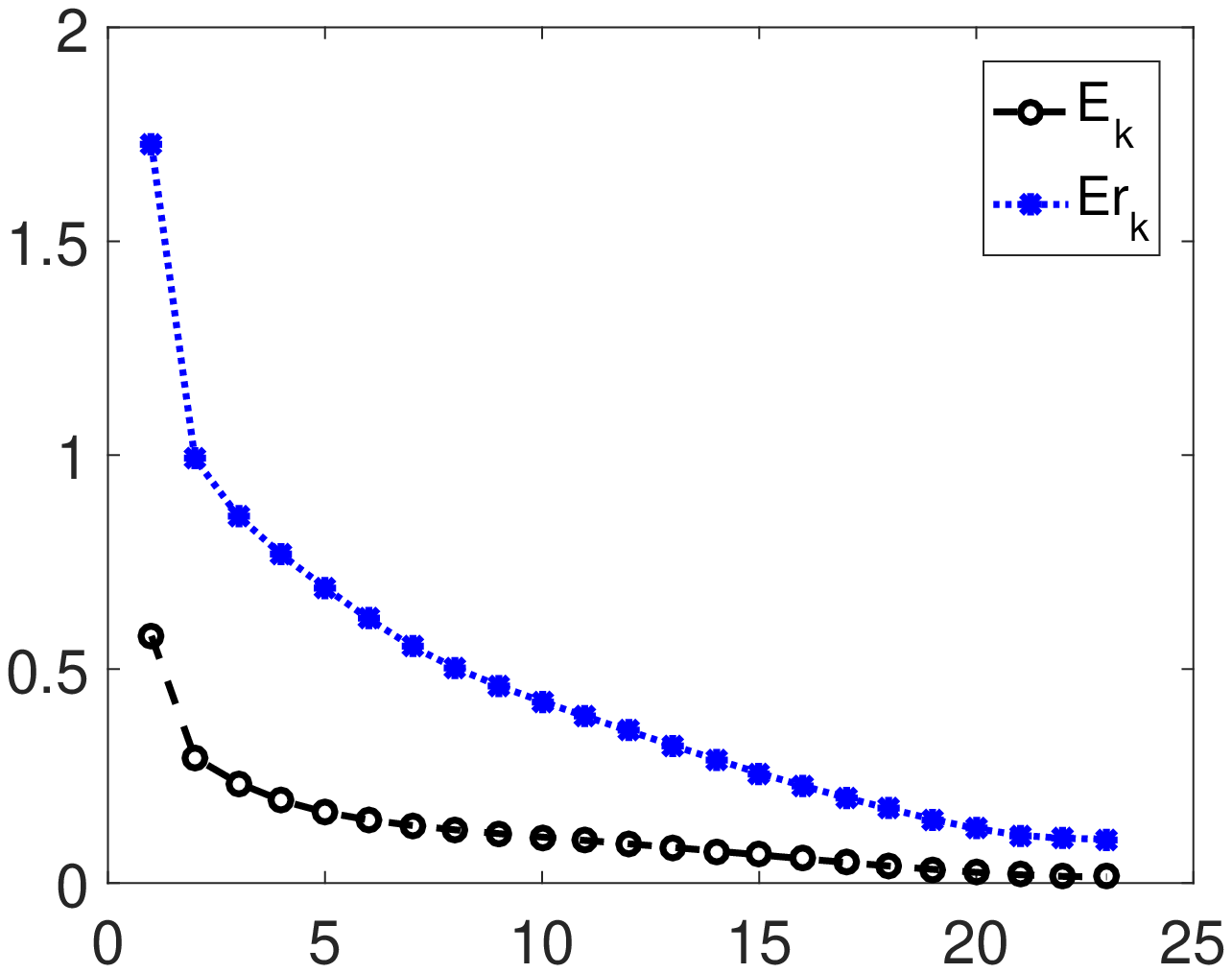}}
	\caption{Reconstructions of a peanut-shaped domain with $1\%$ noise and $\epsilon=0.015$.}\label{4.6}
\end{figure}

\begin{figure}
	\centering 
	\subfigure[5th iteration]{\includegraphics[width=0.32\textwidth]{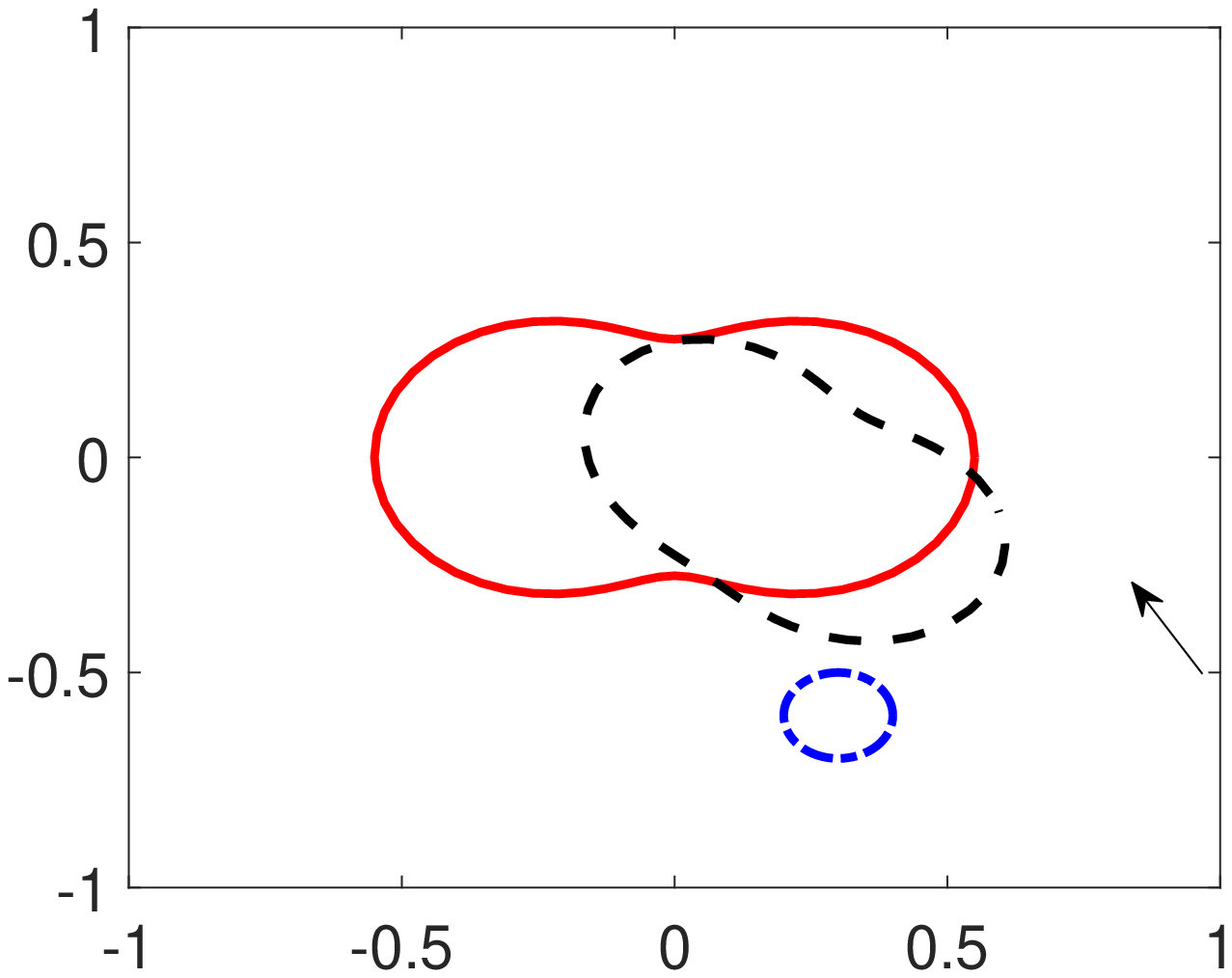}}
	\subfigure[10th iteration]{\includegraphics[width=0.32\textwidth]{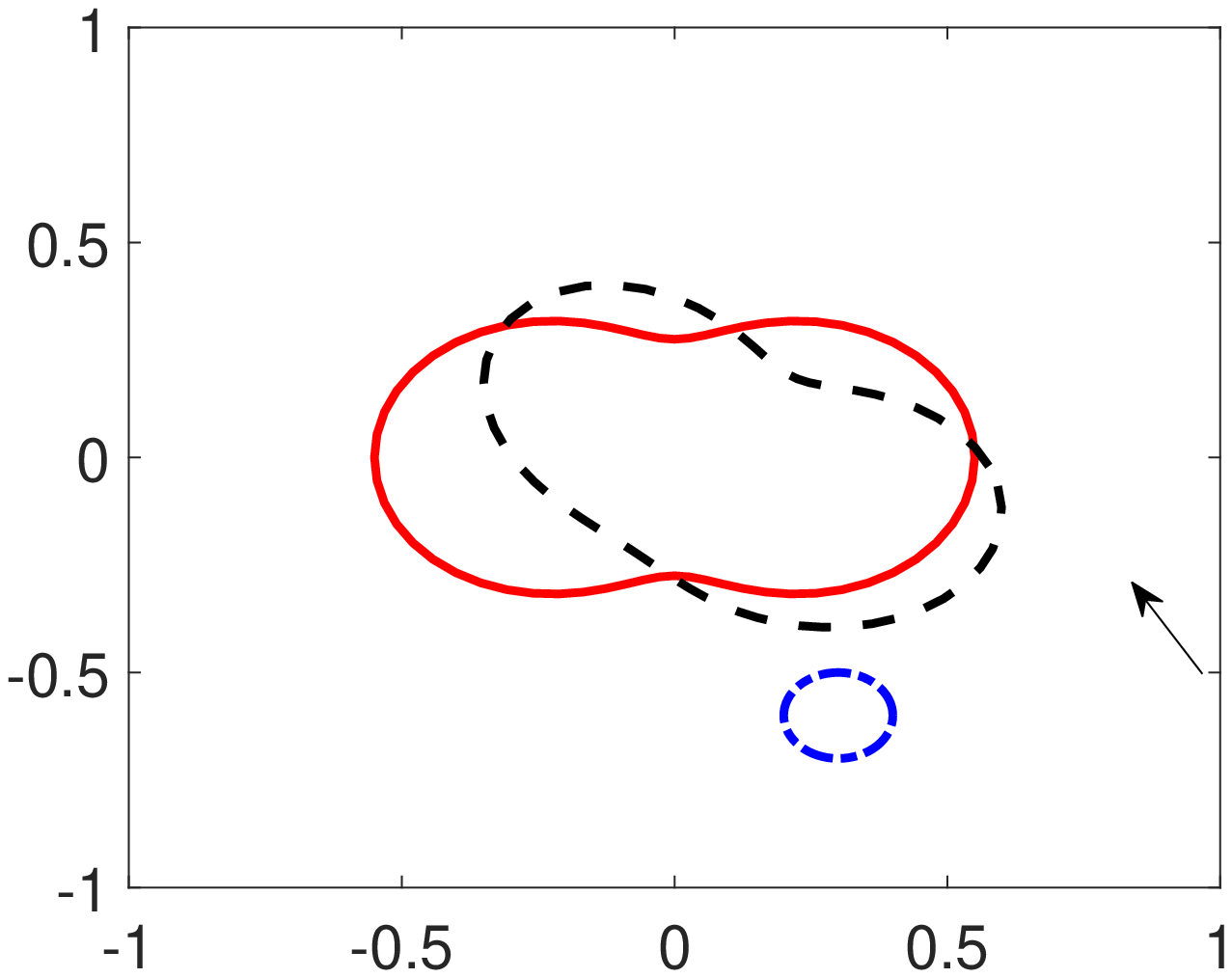}}
	\subfigure[15th iteration]{\includegraphics[width=0.32\textwidth]{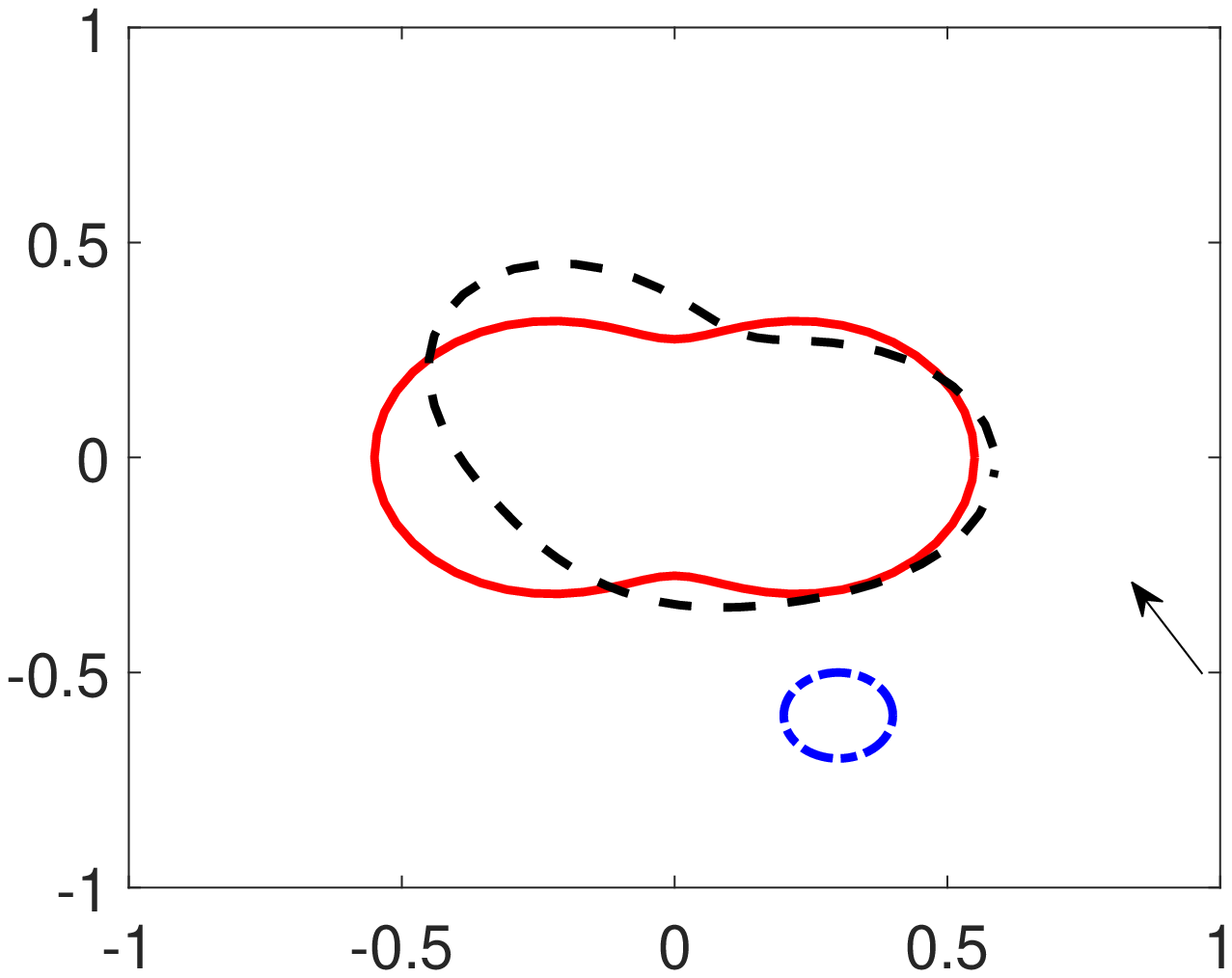}}
	\subfigure[19th iteration]{\includegraphics[width=0.32\textwidth]{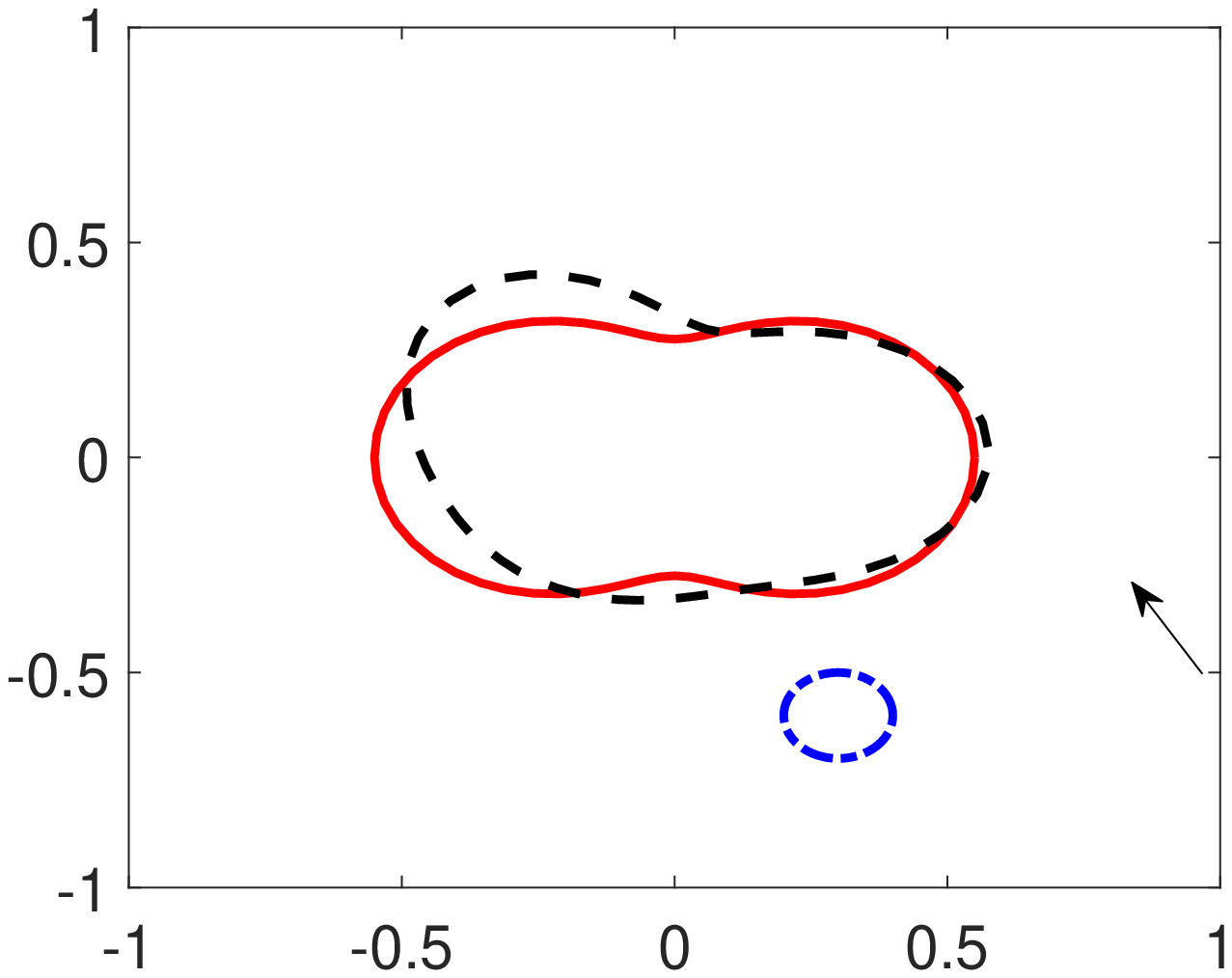}}
	\subfigure[relative error]{\includegraphics[width=0.32\textwidth]{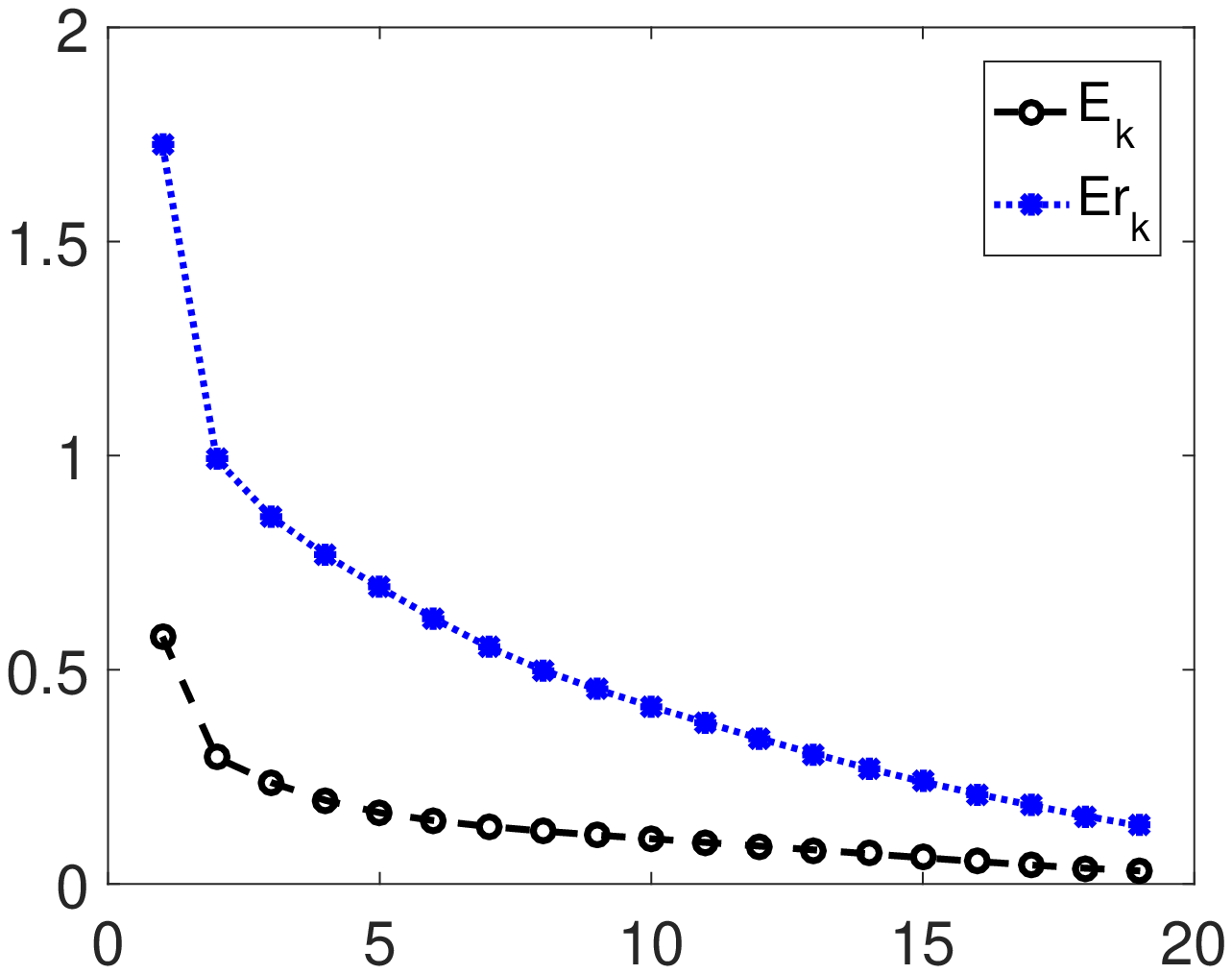}}
	\caption{Reconstructions of a peanut-shaped domain with $5\%$ noise and $\epsilon=0.035$.}\label{4.7}
\end{figure}

\begin{figure}
	\centering 
	\subfigure[]{\includegraphics[width=0.32\textwidth]{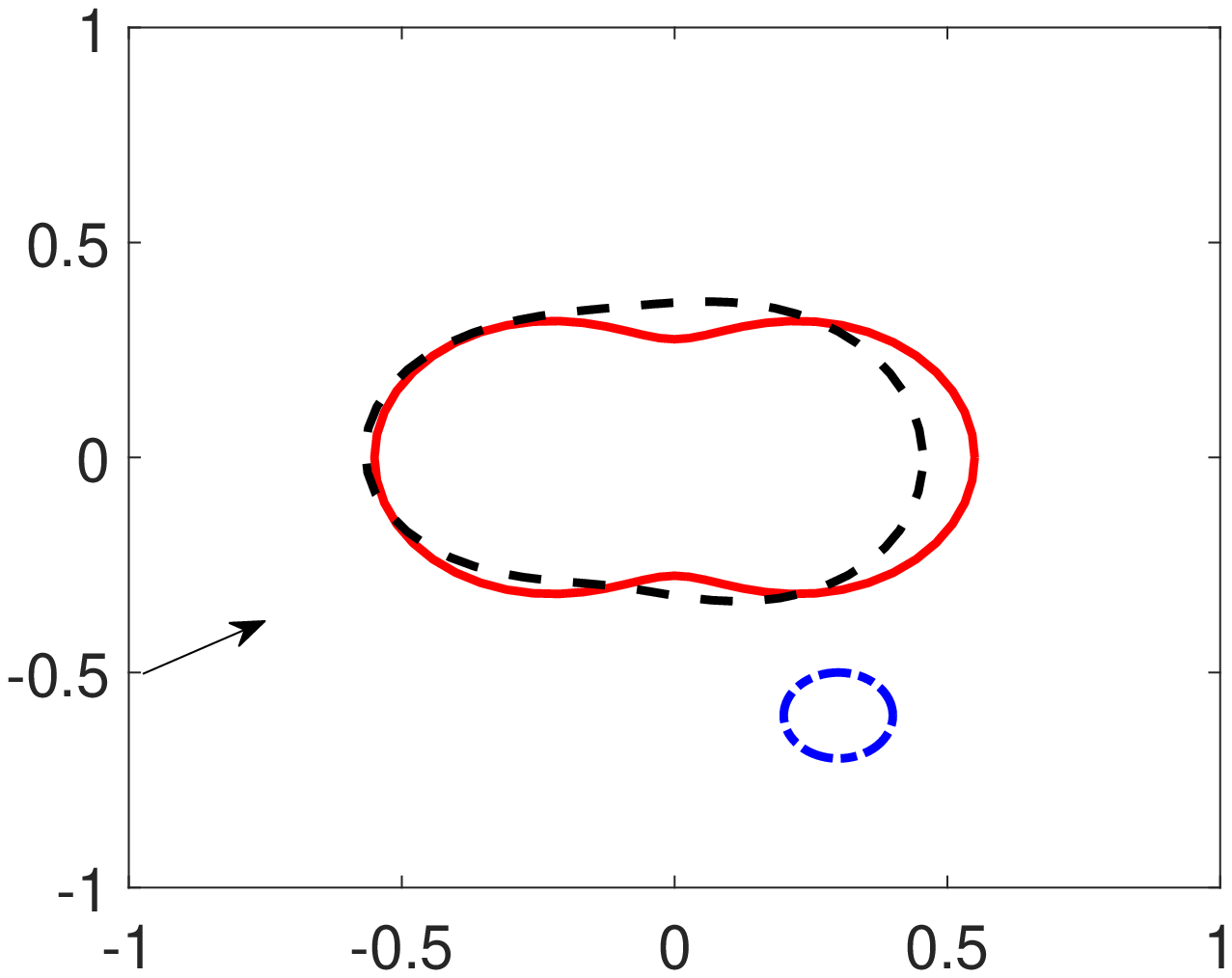}} 
	\subfigure[]{\includegraphics[width=0.32\textwidth]{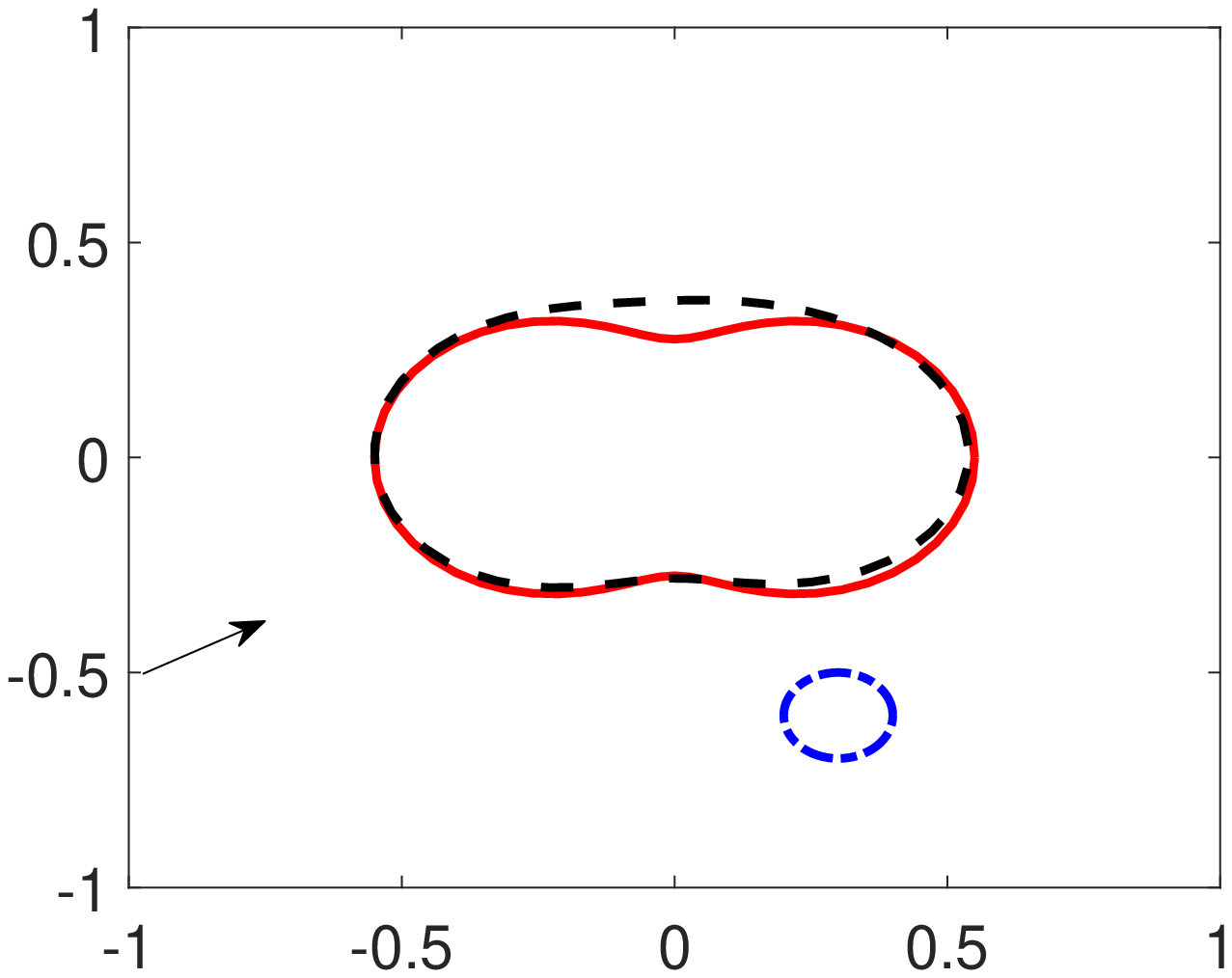}}
	\subfigure[]{\includegraphics[width=0.32\textwidth]{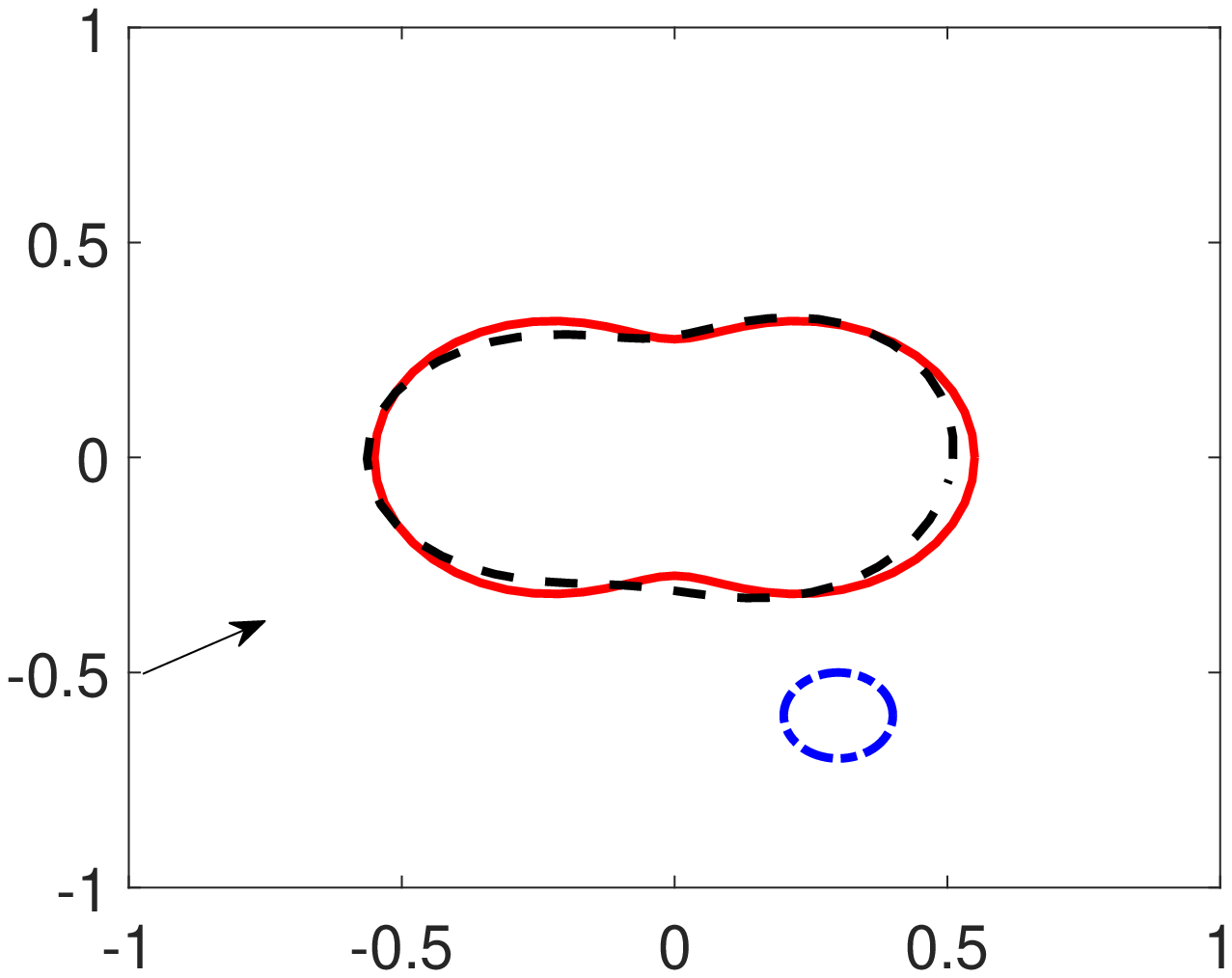}}
	\caption{Reconstructions of a peanut-shaped domain with different reference balls, $1\%$ noise is added. Here, the initial guess $(c_1^{(0)},c_2^{(0)})=(0.3, -0.6), r^{(0)}=0.1$ and $\epsilon=0.015$. (a) $(b_1, b_2)=(3, 0), R=0.2$. (b) $(b_1, b_2)=(4, 0), R=0.4$. (c) $(b_1, b_2)=(5.5, 0), R=0.8$.}\label{4.8}
\end{figure}

In this subsection, we present three numerical examples to illustrate the feasibility of the iterative reconstruction method. To avoid committing an inverse crime, the synthetic data is numerically generated at 64 points, i.e., $n=32$, by the combined double- and single-layer potential method \cite{DR-shu2}. The noisy data
$|u^{\infty,\delta}|^2$ is constructed in the following way
\begin{align*}
|u^{\infty,\delta}|^2=|u^\infty|^2(1+\delta\eta)
\end{align*}
where $\eta$ are normally distributed random numbers ranging in $[-1,1]$ and $\delta>0$ is the relative noise level.

In the numerical examples, we obtain the update $\xi$ from a scaled Newton step with Tikhonov regularization and $H^2$ penalty term, that is
$$
\xi=\rho(\lambda\widetilde{I}+\widetilde{B}^*\widetilde{B})^{-1}\widetilde{B}^*\widetilde{f},
$$
where the scaling factor $\rho\geq0$ is fixed throughout the iterations. According to \cite{OT2007}, the regularization parameters $\lambda$ in equation \eqref{EqualRLHuygens3} are chosen as
\[
\lambda_k:=\Big\||w^\infty|^2-|A^\infty_1(p^{(k-1)}_1,\psi^{(k-1)}_1)
+A^\infty_2(p_2,\psi^{(k-1)}_2)|^2\Big\|_{L^2},\ k=1,2,\cdots.
\]
Note that in each step of the iteration, the derivative $\mathrm{d}r/\mathrm{d}\tau$ is calculated by resorting to the approximation \eqref{updataq} and $r^{(k)}(\tau)=r^{(k-1)}(\tau)+{\Delta r}^{(k-1)}(\tau)$, $k=1,2\cdots$.

Analogously to \cite{RW1997}, the initial approximation is chosen as a circle and the values of the Fourier modes in the directions $\cos\theta$ and $\sin\theta$ (i.e., coefficients $\alpha_1$ and $\beta_1$) are set to be the exact values in the reconstructions. The reason for doing so is to ease the comparison with the exact curve. So the update procedure does not take these two modes into account.

In the subsequent figures, the exact boundary curves are displayed as solid lines, the reconstruction are depicted with the dashed lines $--$, the initial guess are taken to be circles with radius $r^{(0)}=0.1$ indicated by the dash-dotted lines $\cdot-$. The incident directions are denoted by arrows. Throughout all the numerical examples, we set the wavenumber $\kappa=2$, the scaling factor $\rho=0.6$, and the parameter $M=5$.

\begin{figure}
	\centering
	\subfigure[$(c_1^{(0)},c_2^{(0)})=(-0.3, -0.6)$]{\includegraphics[width=0.45\textwidth]{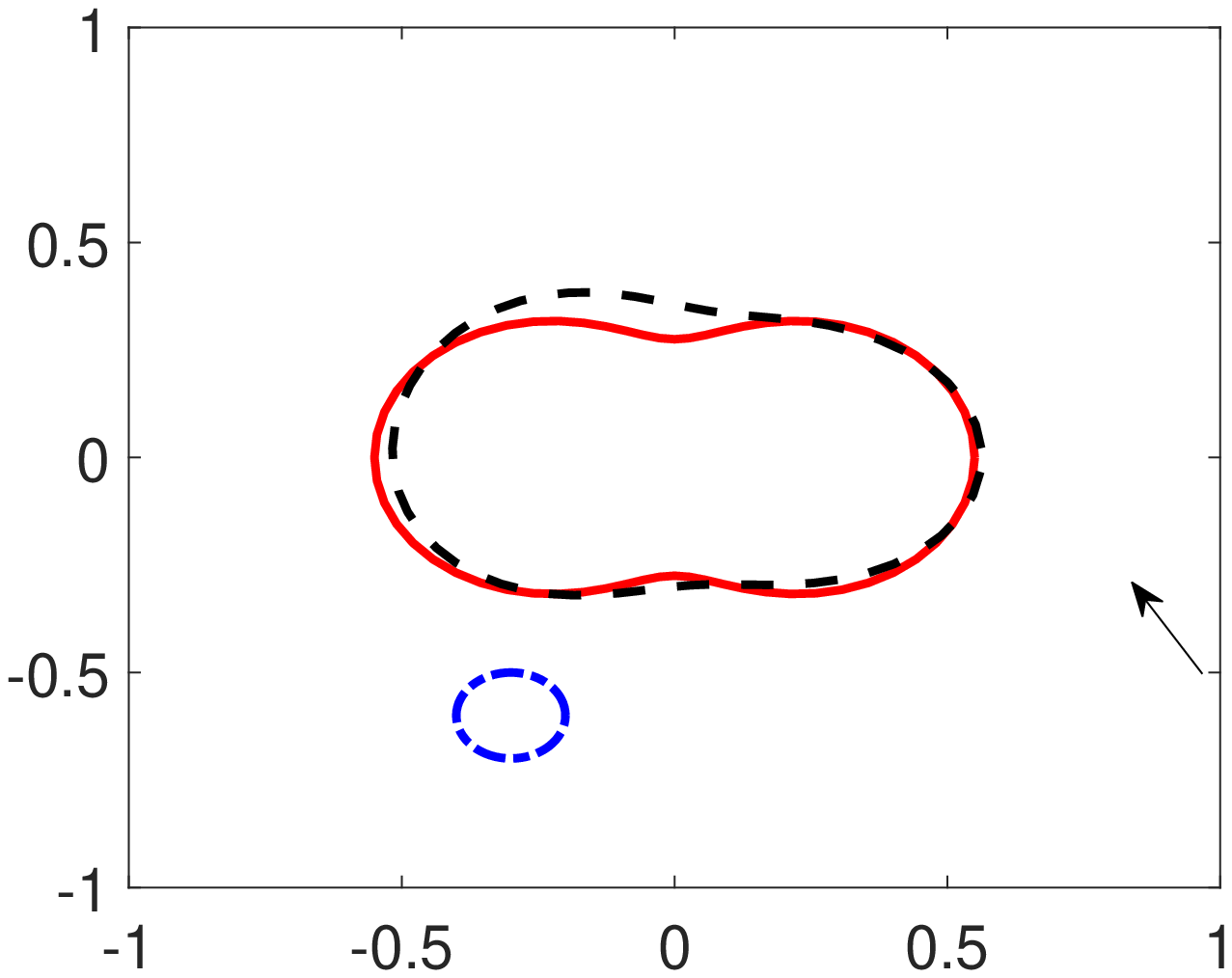}}
	\subfigure[$(c_1^{(0)},c_2^{(0)})=(0.2, 0.7)$]{\includegraphics[width=0.45\textwidth]{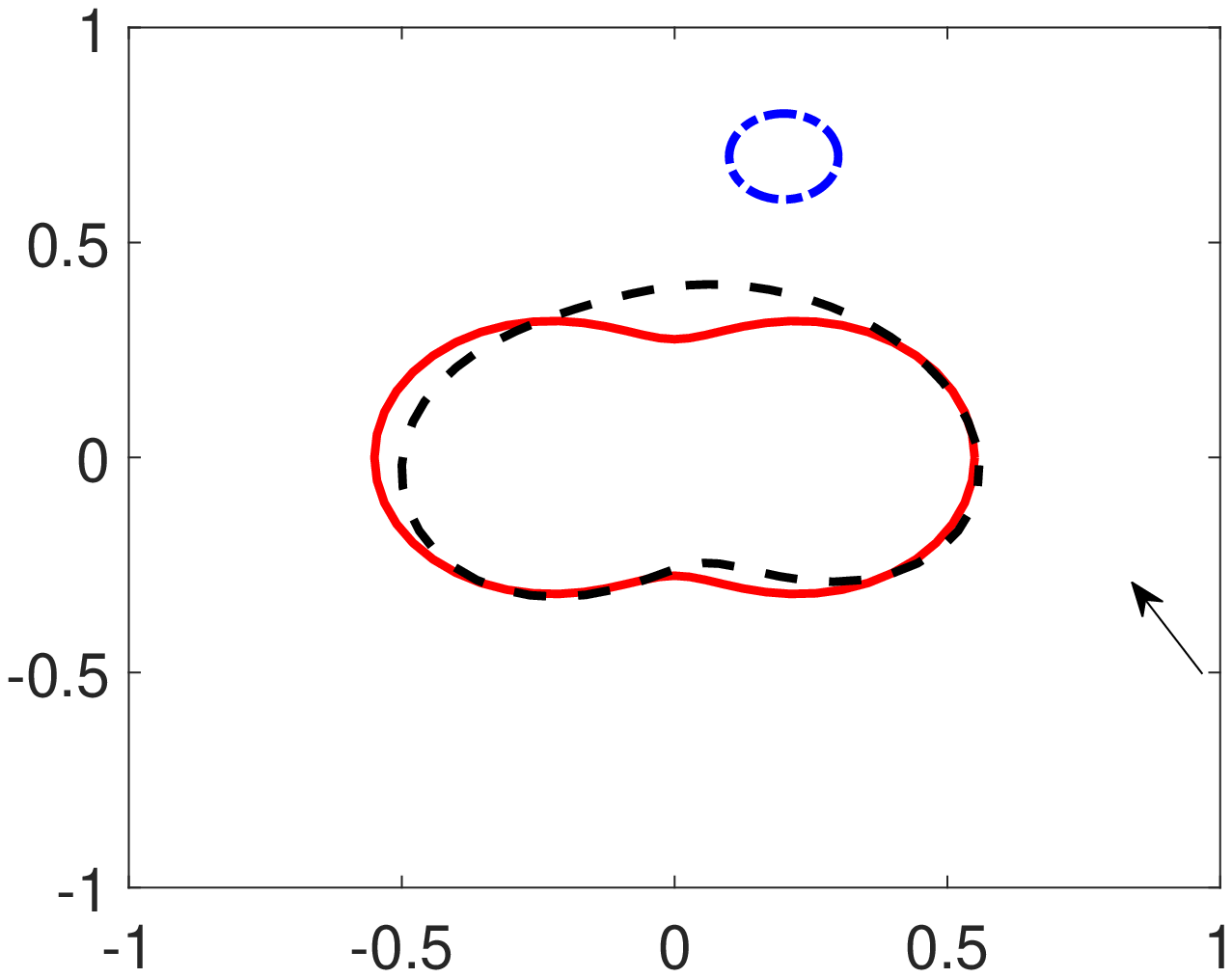}}
	\caption{Reconstructions of a peanut shaped domain with different initial guesses and $1\%$ noise. Here, we are using the reference ball $(b_1, b_2)=(4, 0), R=0.4$ and $\epsilon=0.015$.}\label{4.9}
\end{figure}

\begin{figure}
	\centering 
	\subfigure[$1\%$ noise, $\epsilon=0.015$]{\includegraphics[width=0.45\textwidth]{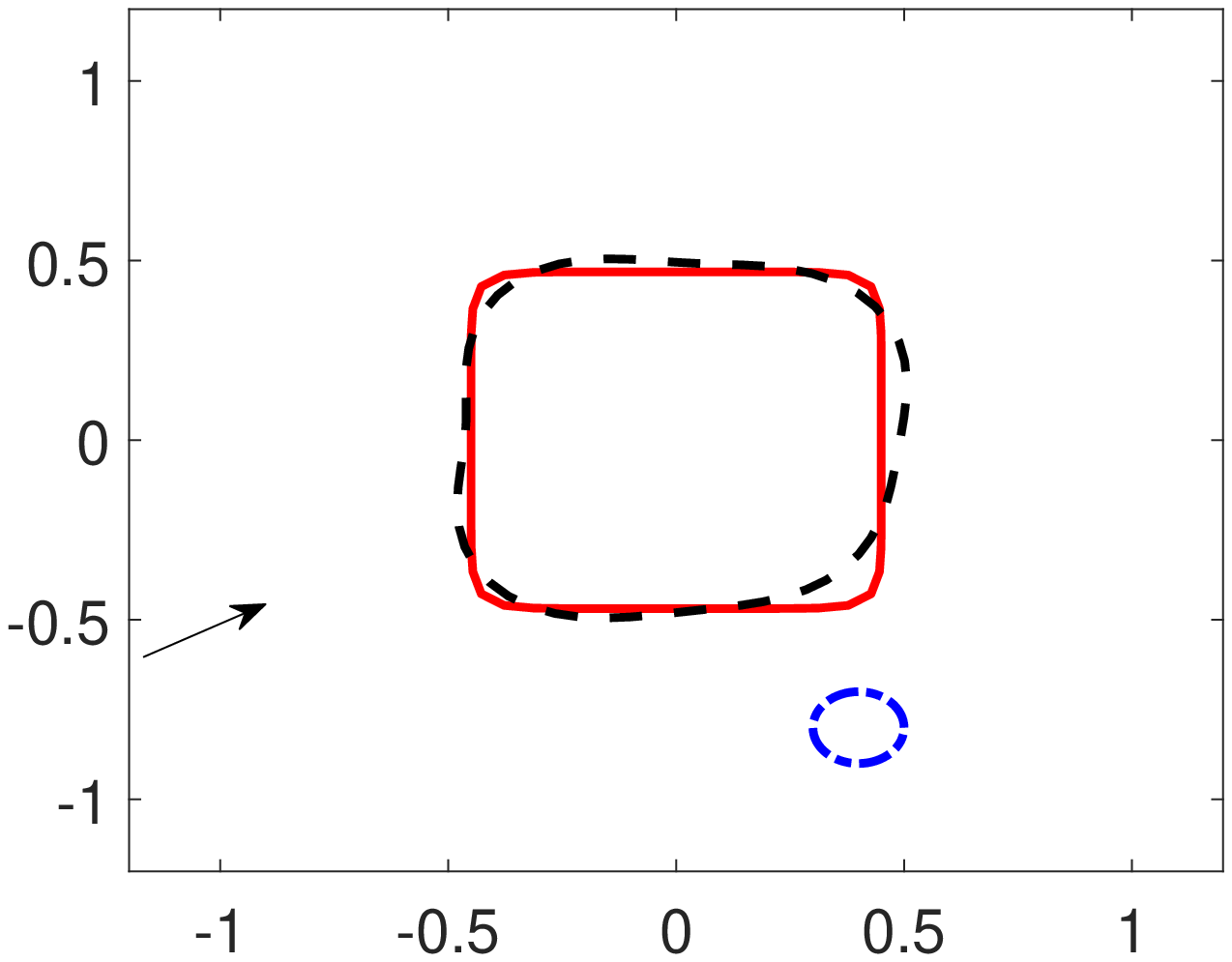}}
	\subfigure[relative error, $1\%$ noise]{\includegraphics[width=0.45\textwidth]{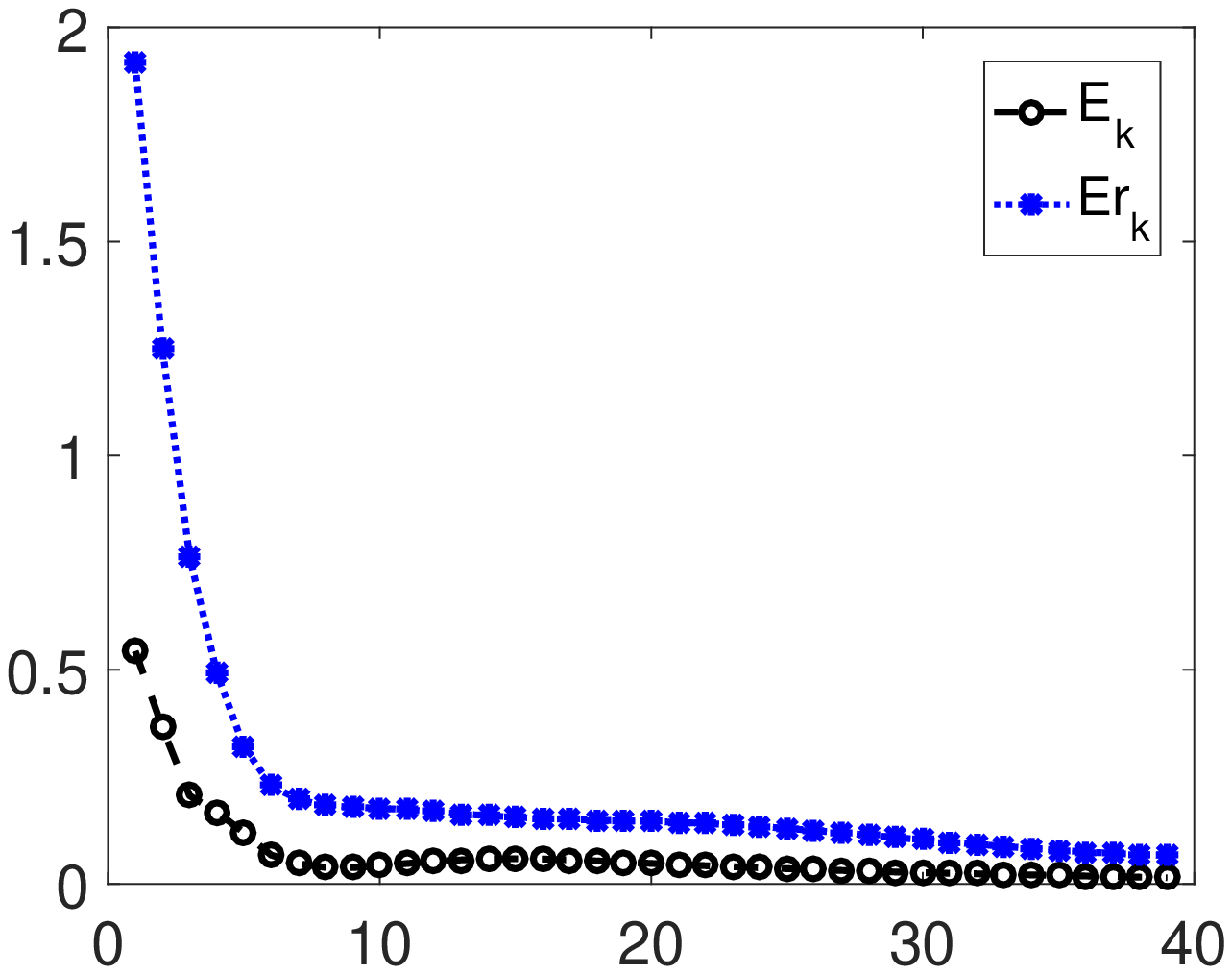}}
	\subfigure[$5\%$ noise, $\epsilon=0.035$]{\includegraphics[width=0.45\textwidth]{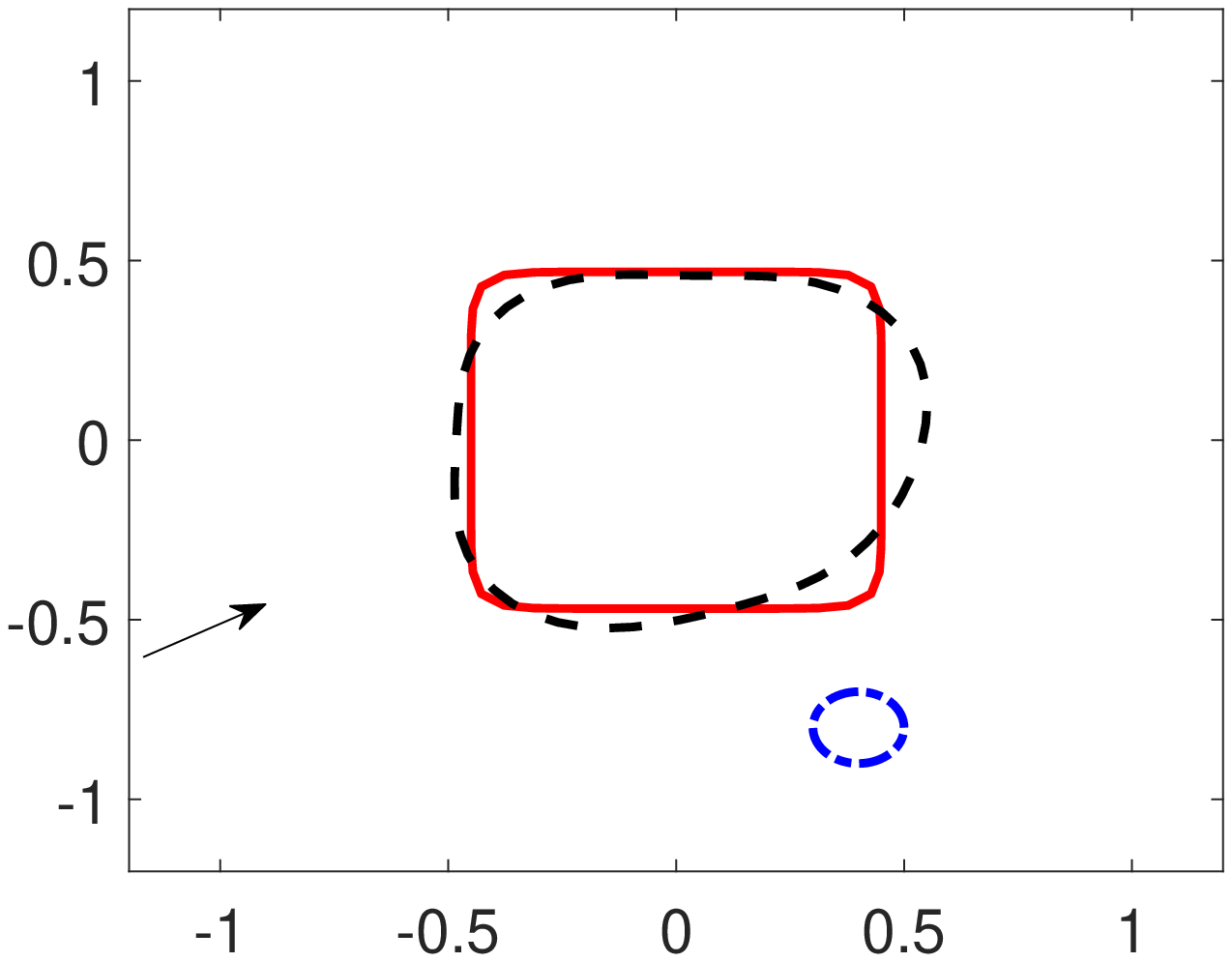}}
	\subfigure[relative error, $5\%$ noise]{\includegraphics[width=0.45\textwidth]{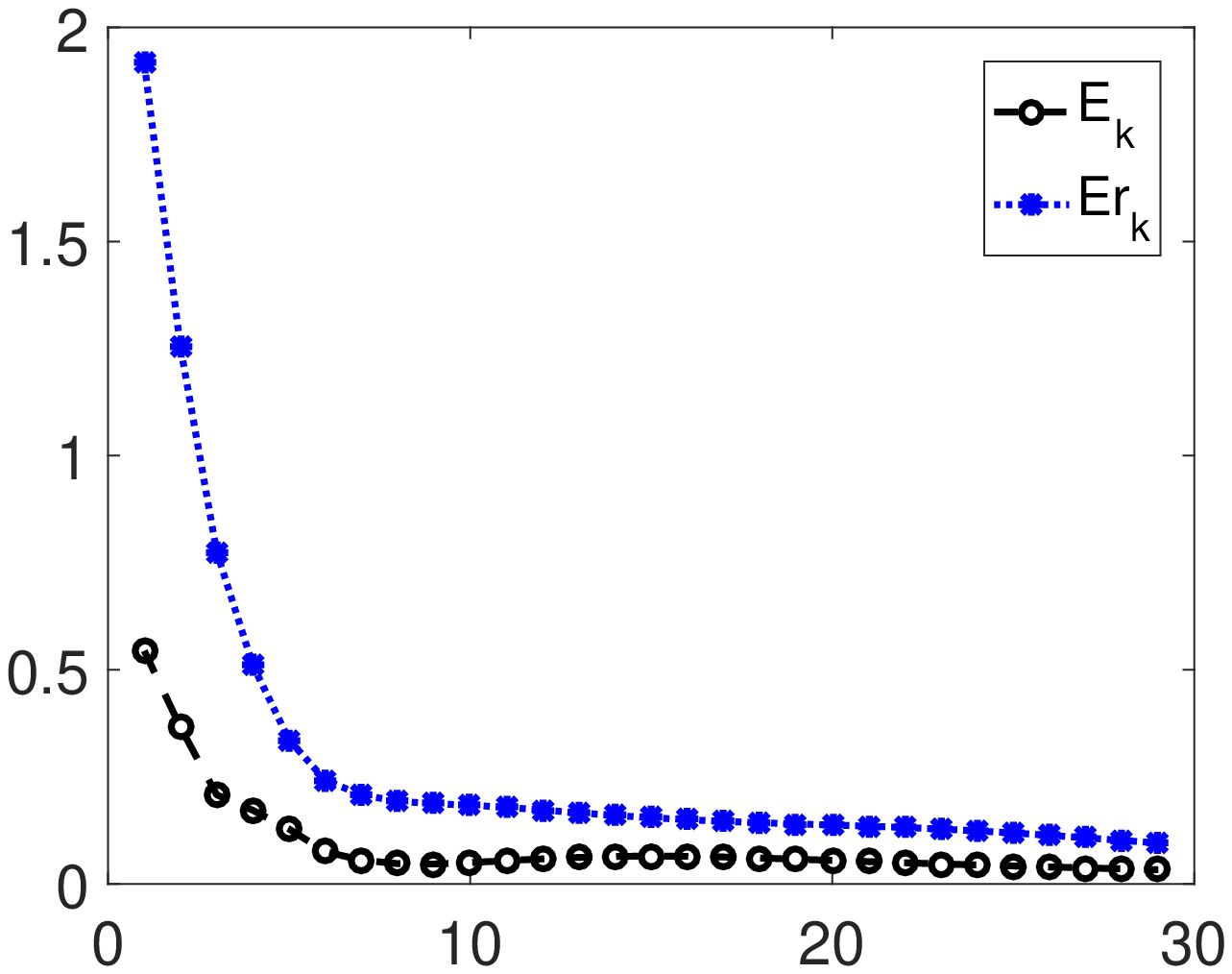}}
	\caption{Reconstructions of a rectangle shaped domain with $1\%$ and $5\%$ noise data respectively. Here, we are using the initial guess $(c_1^{(0)},c_2^{(0)})=(0.4, -0.8), r^{(0)}=0.1$ and the reference ball $(b_1, b_2)=(4, 0), R=0.5$.}\label{4.10}
\end{figure}

\begin{figure}
	\centering 
	\subfigure[$(c_1^{(0)}, c_2^{(0)})=(0.4, -0.8)$]{\includegraphics[width=0.45\textwidth]{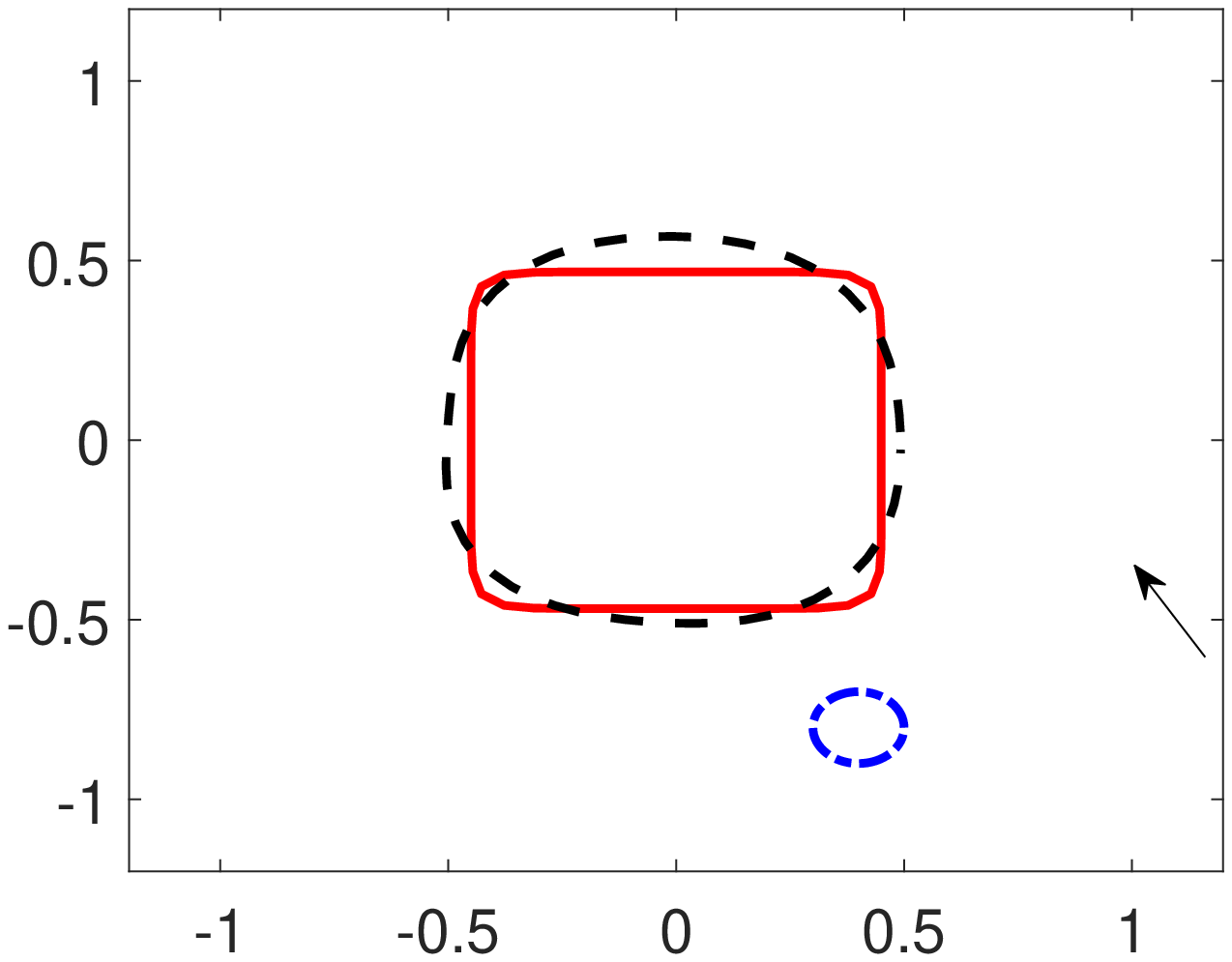}}
	\subfigure[$(c_1^{(0)}, c_2^{(0)})=(-0.3, 0.7)$]{\includegraphics[width=0.45\textwidth]{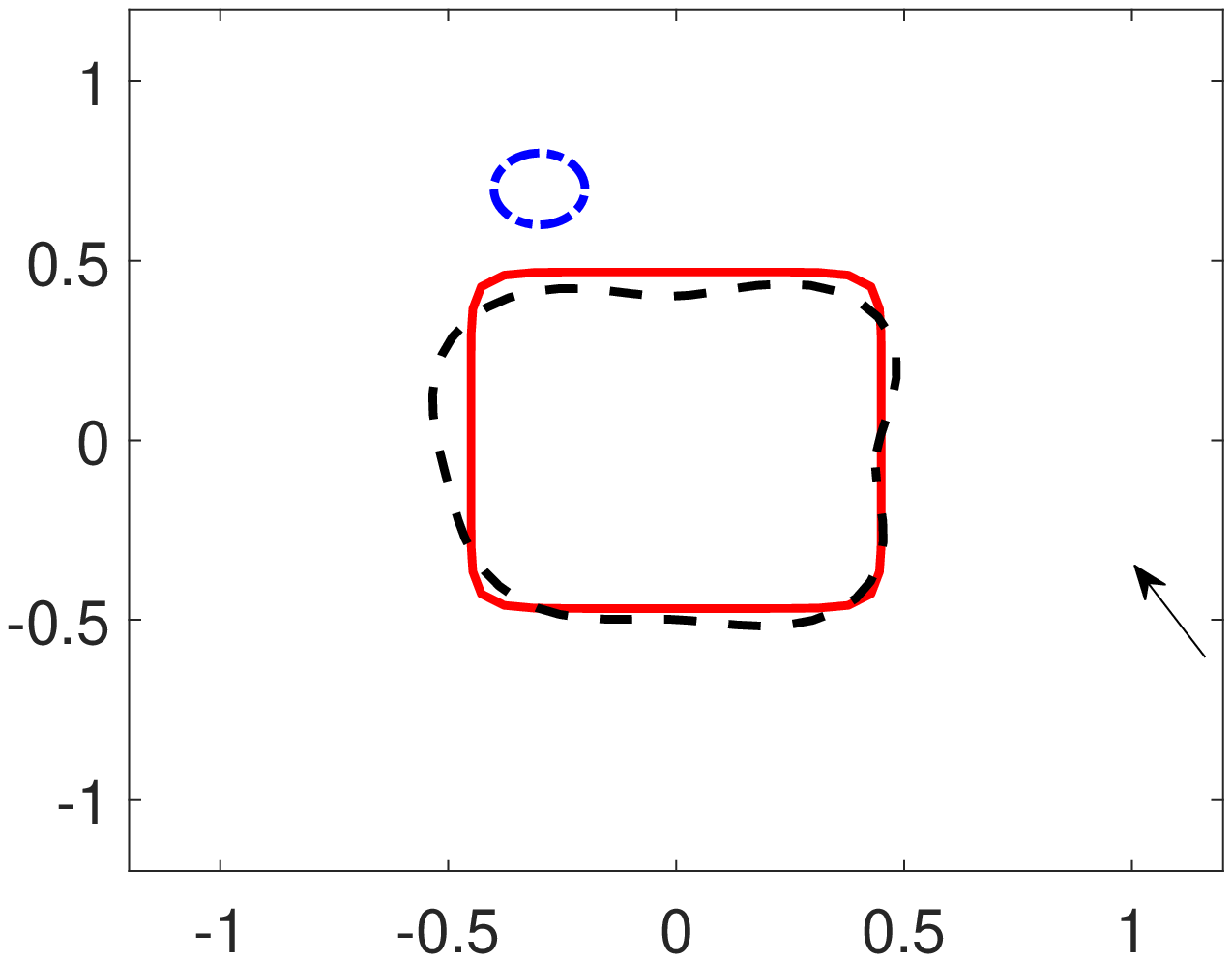}}
	\caption{Reconstructions of a rectangle shaped domain with different initial guesses ($1\%$ noise is added). Here, we are using the incoming wave direction $d=(\cos(2\pi/3),\sin(2\pi/3))$, the reference ball $(b_1, b_2)=(4, 0), R=0.5$ and $\epsilon=0.015$.}\label{4.11}
\end{figure}

\begin{figure}
	\centering 
	\subfigure[ $(b_1, b_2)=(4, 0), R=0.5$]{\includegraphics[width=0.45\textwidth]{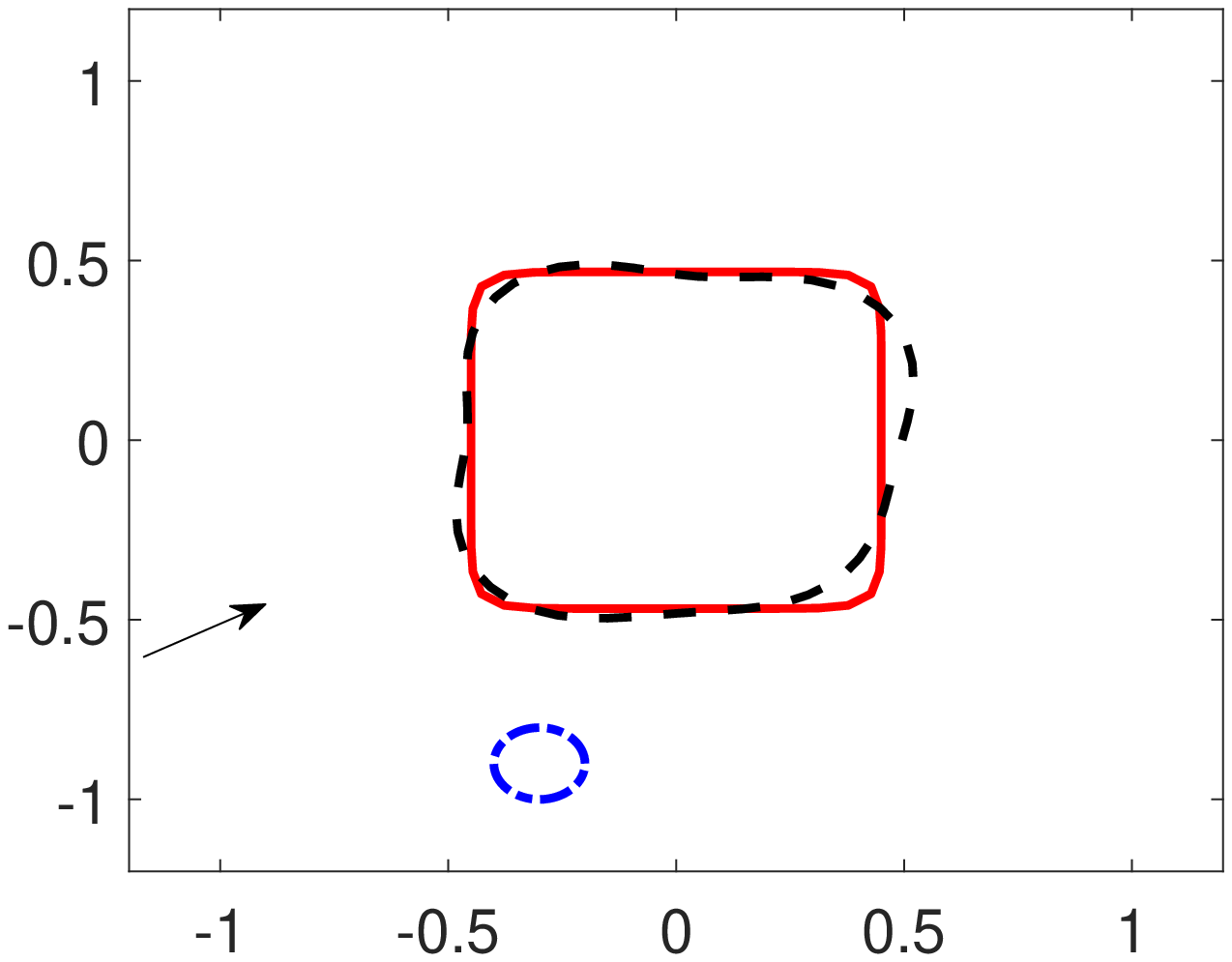}} 
	\subfigure[$(b_1, b_2)=(7, 0), R=0.8$]{\includegraphics[width=0.45\textwidth]{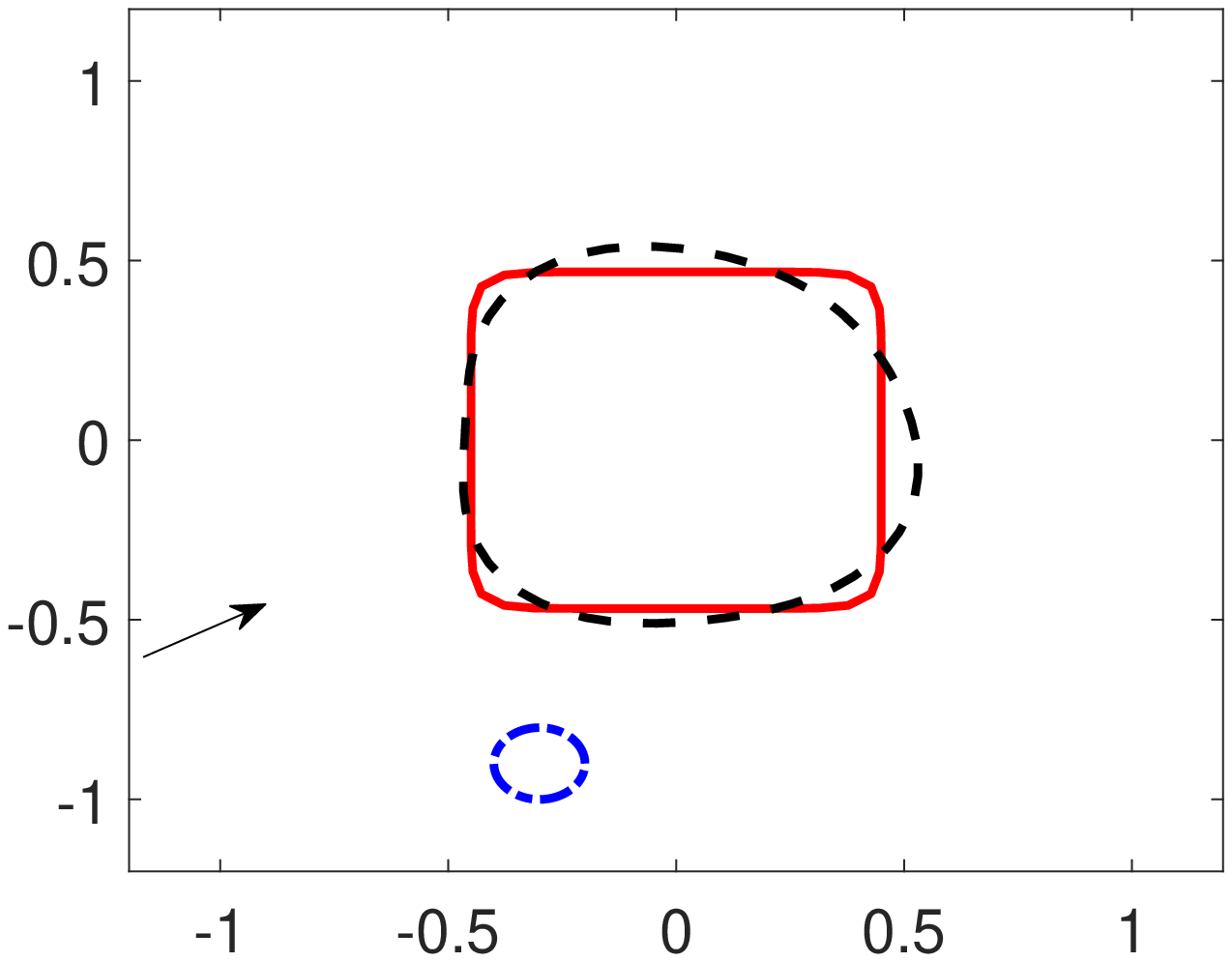}}
	\caption{Reconstructions of a rectangle-shaped domain with different reference balls, $1\%$ noise is added. Here, the incoming wave direction $d=(\cos(\pi/6), \sin(\pi/6))$, the initial guess $(c_1^{(0)}, c_2^{(0)})=(-0.3, -0.9), r^{(0)}=0.1$, and $\epsilon=0.015$.}\label{4.12}
\end{figure}

\begin{example}\label{E1}
	\rm In the first example, we aim to reconstruct an apple-shaped obstacle with parametrized boundary
	$$
	p_1(t)=\frac{0.55(1+0.9\cos{t}+0.1\sin{2t})}{1+0.75\cos{t}}(\cos{t},\sin{t}), \quad t\in [0,2\pi].
	$$
	
	Here the incident direction $d=(\cos(-\pi/6),\sin(-\pi/6))$ is used in Fig. \ref{4.1}-\ref{4.3}. Several snapshots of the iterative process are shown in Fig. \ref{4.1} and Fig. \ref{4.2}, where the initial guess is the circle centered at $(-0.7, 0.45)$ with radius 0.1, and the reference ball is centered at $(4,0)$ with radius 0.4. Moreover, the relative $L^2$ error $Er_k$ between the reconstructed and exact boundaries and the error $E_k$ defined in \eqref{relativeerror} are also presented with respect to the number of iterations. As we can see from the
	figures, the trend of two error curves is basically the same for larger number of iteration. So, the choice of the stopping criteria is reasonable.  The reconstructions with different reference balls for the incoming wave direction $d=(\cos(-\pi/6),\sin(-\pi/6))$ and $d=(\cos(4\pi/3),\sin(4\pi/3))$ are shown in Fig. \ref{4.3} and Fig. \ref{4.5} respectively, and the reconstructions with different incoming wave directions are presented in Fig. \ref{4.4}. As shown in these results, the location and shape of the obstacle could be simultaneously and satisfactorily reconstructed.
\end{example}

\begin{example}\label{E2}
	\rm For the second example, we consider the scattering by a peanut-shaped obstacle described by
	$$
	p_1(t)=0.275\sqrt{3\cos^2{t}+1}(\cos{t},\sin{t}), \quad
	t\in[0,2\pi].
	$$
	
	In this example, the reconstructions with $1\%$ noise and $5\%$ noise, corresponding to the incoming wave direction $d=(\cos(2\pi/3),\sin(2\pi/3))$, are shown in Fig. \ref{4.6} and
	Fig. \ref{4.7}, respectively. Here we use the initial guess $(c_1^{(0)},c_2^{(0)})=(0.3, -0.6), r^{(0)}=0.1$ and the reference ball $(b_1, b_2)=(4, 0), R=0.4$. The relative $L^2$ error $Er_k$ and the error $E_k$ are also presented in the figures. The reconstructions with different reference balls for the incoming wave direction $d=(\cos(\pi/6),\sin(\pi/6))$ are shown in Fig. \ref{4.8}, and the reconstructions with different initial guesses for  the incoming wave direction $d=(\cos(2\pi/3),\sin(2\pi/3))$ are shown in Fig. \ref{4.9}.
\end{example}

\begin{example}\label{E3}
	\rm Our third example is intended to reconstruct a rounded rectangle
	given by
	$$
	p_1(t)=\frac{9}{20}\left(\cos^{10}t+\frac{2}{3}\sin^{10}t\right)^{-1/10}(\cos{t},\sin{t}),
	\quad t\in[0,2\pi].
	$$
	
	In this example, the reconstructed obstacle and the relative error $Er_k$ and $E_k$ from the phaseless far-field data with $1\%$ noise and $5\%$ noise, corresponding to the incoming wave direction $d=(\cos(\pi/6),\sin(\pi/6))$, are shown in Fig. \ref{4.10}. The relative $L^2$ error $Er_k$ and the error $E_k$ are also presented in the figures. The influences of the choices of initial guesses and reference balls are shown in Fig. \ref{4.11} and Fig. \ref{4.12} respectively.
\end{example}

The above numerical results illustrate that by adding a reference ball artificially to the inverse scattering system the iteration method gives a feasible reconstruction of the location and the shape of the obstacle from phaseless far-field data for one incident field. The reference ball causes few extra computational costs, but breaks the translation invariance and brings information about the location of the obstacle. In addition, a promising feature of the algorithm is that the Fr\'{e}chet derivatives involved in our method can be explicitly characterized as integral operators and thus easily evaluated. Hence it does not require the solution of the forward problem in each iteration step and it is very easy to implement with computational efficiency. To evaluate the computational time, our codes of the experiments are written in Matlab and run on a laptop with 2.6 GHz CPU. According to the fact that the CPU time for all the reconstructions is less then 10 seconds, we conclude that our algorithm is very fast. 


\section {Conclusions and future works}

In this paper, a new numerical method is devised to solve the inverse obstacle scattering problem from the modulus of the
far-field data for one incident field. That is, by introducing a reference ball artificially to the inverse scattering system, the translation invariance of the phaseless far-field pattern can be broken down, and an iterative scheme which is based on a system of boundary integral equations is propose to reconstruct both the location and shape of the obstacle. The reference ball causes few extra computational costs, but breaks the translation invariance and brings the location information of the obstacle. The numerical implementation details of the iterative scheme is described, and the numerical examples illustrate that the iterative method yields satisfactory reconstructions. 

Concerning our future work, the proposed methodology could be extended directly to the case of recovering a sound-hard or impedance obstacle, as well as the three-dimensional case. Although our numerical results show that the proposed novel iterative scheme works very well, the corresponding theoretical justifications of convergence and stability are still open. In other words, the mathematical analysis of this method is beyond the scope of our current work and deserves future investigations. In addition, we would also like to study the applicability of this approach to imaging crack-like scatterers from phaseless data. Moreover, we believe that the reference ball based iteration approach is also a feasible technique for solving the phaseless inverse electromagnetic scattering problem.

\section*{Acknowledgements}

The first author is supported by NSFC grant 11771180. The second author is supported by NSFC grant 11671170. The third author is supported by NSFC grants 11601107, 41474102 and 11671111.


\end{document}